\documentclass[12pt]{amsart}

\usepackage[utf8]{inputenc}
\usepackage[T1]{fontenc}
\usepackage{lmodern}
\usepackage{amsmath,amssymb,amsthm,mathtools}
\usepackage[a4paper,margin=2.7cm]{geometry}
\usepackage{xcolor}
\usepackage{tikz-cd}
\usepackage{booktabs}
\usepackage{enumitem}
\usepackage[colorlinks=true,linkcolor=black,citecolor=black,urlcolor=black]{hyperref}

\newtheorem{theorem}{Theorem}[section]
\newtheorem{proposition}[theorem]{Proposition}
\newtheorem{lemma}[theorem]{Lemma}
\newtheorem{corollary}[theorem]{Corollary}
\theoremstyle{definition}
\newtheorem{definition}[theorem]{Definition}
\newtheorem{example}[theorem]{Example}
\theoremstyle{remark}
\newtheorem{remark}[theorem]{Remark}

\numberwithin{equation}{section}
\newcommand{\R}{\mathbb{R}}
\newcommand{\g}{\mathfrak{g}}
\newcommand{\X}{\mathfrak{X}}
\newcommand{\Hor}{\mathcal{H}}
\newcommand{\Ver}{\mathcal{V}}
\newcommand{\Reeb}{\mathcal{R}}
\newcommand{\Lie}{\mathcal{L}}
\newcommand{\pd}[2]{\frac{\partial #1}{\partial #2}}
\newcommand{\Ad}{\mathrm{Ad}}
\newcommand{\ad}{\mathrm{ad}}
\newcommand{\id}{\mathrm{id}}
\newcommand{\sharpL}{\sharp_\Lambda}
\newcommand{\perpd}{\perp_{d\eta}}
\newcommand{\perpL}{\perp_\Lambda}
\renewcommand{\Im}{\operatorname{Im}}

\makeatletter
\def\@settitle{\begin{center}\baselineskip16\p@\relax\bfseries\Large\@title\end{center}}
\makeatother

\title[Hamilton--Jacobi theory and reduction for contact systems]{Hamilton--Jacobi theory, symmetries and coisotropic reduction for contact Hamiltonian systems}


\begin{document}

\renewcommand{\shortauthors}{M. DE LE\'ON AND M. E. EYREA IRAZ\'U}

\maketitle

\begin{center}
\small
\textsc{Manuel de Le\'on}\\[2pt]
Instituto de Ciencias Matem\'aticas (CSIC--UAM--UC3M--UCM),\\ C/ Nicol\'as Cabrera 13--15, Campus Cantoblanco, 28049 Madrid, Spain\\
Real Academia de Ciencias Exactas, F\'isicas y Naturales de Espa\~na,\\ C/ Valverde 22, 28004 Madrid, Spain\\
Great Bay University, Dongguan, Guangdong, China\\
\texttt{mdeleon@icmat.es}\\[8pt]
\textsc{Mar\'ia Emma Eyrea Iraz\'u}\\[2pt]
Departamento de Matem\'atica, Centro de Matem\'atica de La Plata (CMaLP),\\ Facultad de Ciencias Exactas, Universidad Nacional de La Plata, La Plata, Argentina\\
\texttt{maemma@mate.unlp.edu.ar}
\end{center}

\begin{abstract}
We study how the Hamilton--Jacobi theory for contact Hamiltonian systems interacts with symmetries and coisotropic reduction, extending to the contact setting the programme carried out by de Le\'on, Mart\'in de Diego and Vaquero for symplectic Hamiltonian systems. We show that a section of $T^*Q\times\R\to Q\times\R$ has coisotropic image if and only if its image is foliated by Legendrian submanifolds, locally a one-parameter family of $1$-jets; accordingly, the Hamilton--Jacobi equation for the evolution vector field says that $H$ is constant on each leaf, and the one for the Hamiltonian vector field becomes the evolutionary equation $\partial_tS+H(q,\partial_qS,S)=0$. For actions by contactomorphisms, momentum maps are dissipated quantities, so that no Legendrian solution has constant non-zero momentum when $\Reeb(H)$ vanishes nowhere, and a closed connected Legendrian submanifold is $G$-invariant if and only if it lies in $J^{-1}(0)$. Using contact coisotropic reduction, we establish a one-to-one correspondence between $G$-invariant solutions for $H$ and solutions for the reduced Hamiltonian on $T^*(Q/G)\times\R$. Non-zero momenta $\mu$ vanishing on $[\g,\g]$ are handled in two ways: by a twisted reduction, valid when $\Reeb(H)=0$ on $J^{-1}(\mu)$ and conjugated to the reduction at zero by a strict contactomorphism, and, for arbitrary dissipative Hamiltonians, by a reduction along the ray $\R^+\mu$ by the subgroup with Lie algebra $\ker\mu$. No magnetic terms appear. Several examples with linear dissipation illustrate the results.
\end{abstract}

\tableofcontents

\section{Introduction}

The Hamilton--Jacobi equation is one of the cornerstones of classical mechanics \cite{AM78,Ar89}. In its time-independent form, $H(q^i,\partial S/\partial q^i)=E$, it allows one to reduce the integration of Hamilton's equations to the solution of a single first order partial differential equation, and a complete solution provides a canonical transformation that trivializes the dynamics. Besides its classical relevance, the equation links mechanics with geometric optics and quantum mechanics, and it is at the origin of the construction of symplectic integrators through generating functions.

During the last two decades a \emph{geometric} Hamilton--Jacobi theory has been developed \cite{CGMMMR06,BLdLMdD13}. Its starting point is the observation that, for a $1$-form $\gamma$ on $Q$, the Hamiltonian vector field $X_H$ on $T^*Q$ and its projection $X_H^\gamma=T\pi_Q\circ X_H\circ\gamma$ are $\gamma$-related if and only if $d(H\circ\gamma)=0$, provided that $\gamma$ is closed, i.e., that $\Im\gamma$ is a Lagrangian submanifold. This point of view, which puts the emphasis on Lagrangian submanifolds rather than on generating functions, has been extended to nonholonomic systems \cite{IdLMdD08,dLMMdD10,CGMMMR10}, Poisson manifolds \cite{dLMdDV14} and many other settings. In \cite{dLMdDV17}, de Le\'on, Mart\'in de Diego and Vaquero studied how this geometric Hamilton--Jacobi theory fits with the reduction theory of Hamiltonian systems with symmetry \cite{MW74,MMOPR07}. Using coisotropic reduction of Lagrangian submanifolds \cite{We77} and cotangent bundle reduction, they obtained a reduction and reconstruction procedure for (generalized) solutions of the Hamilton--Jacobi equation, which unifies and extends previous work by Ge and Marsden \cite{GM88} and recovers the classical ans\"atze for cyclic coordinates; see also \cite{FdLMMdDV17} for the role of symmetries in a groupoid approach to Hamilton--Jacobi theory and generating functions. The reduction of the Hamilton--Jacobi problem has also been studied for nonholonomic systems \cite{EdLJS}, for systems with external forces \cite{dLLLG22} (see also \cite{LG24}), and in the framework of the extended Hamilton--Jacobi theory \cite{GMP21}.

On the other hand, contact Hamiltonian systems have attracted a great deal of attention in recent years \cite{BCT17,Br17,dLL19,dLLrev19,GGMRR20,CCM18}. They provide a geometric framework for mechanical systems with dissipation, whose equations of motion can be derived from the variational principle of Herglotz \cite{He30}, and they provide a natural setting for equilibrium and non-equilibrium thermodynamics \cite{MNSS91,vdSM18,SdLLMdD20} (see \cite{dLB25} for an alternative description in terms of almost cosymplectic structures). A contact Hamiltonian $H$ on a contact manifold $(M,\eta)$ defines two relevant vector fields: the Hamiltonian vector field $X_H$, whose flow consists of conformal contactomorphisms and does not preserve $H$, and the evolution vector field $\mathcal{E}_H=X_H+H\Reeb$, which preserves $H$ and plays a central role in thermodynamics. The Hamilton--Jacobi theory for contact Hamiltonian systems was initiated in \cite{dLS17} (see also \cite{GP20}) and developed systematically in \cite{dLLMB21} (see also the review \cite{EdLLSZ22}, and \cite{CdLEILG24} for hybrid and forced systems), where four Hamilton--Jacobi equations were obtained: two for $X_H$ and two for $\mathcal{E}_H$, depending on whether one considers sections of $T^*Q\times\R\to Q\times\R$ (with coisotropic image) or sections of $T^*Q\times\R\to Q$ (with Legendrian image). See also \cite{dLLLGR23} for the relation with integrability, \cite{GG22} for a Hamilton--Jacobi theory for contact Hamiltonians understood as sections of a line bundle (which also covers non-trivial contact structures), and \cite{CCJWY20} for an analytical approach to the evolutionary equation $\partial_t u+H(x,\partial_x u,u)=0$.

Symmetries of contact Hamiltonian systems behave differently from the symplectic ones. If a Lie group $G$ acts by contactomorphisms leaving $H$ invariant, the components of the momentum map are not conserved but \emph{dissipated} quantities \cite{dLL20,GGMRR20}; compare with the behaviour of momentum maps for mechanical systems with Rayleigh dissipation modelled by external forces \cite{dLLLG21}. Contact reduction has been studied by Albert \cite{Al89}, Geiges \cite{Ge97}, Loose \cite{Lo01}, Willett \cite{Wi02} and Zambon and Zhu \cite{ZZ06}, cosphere bundle reduction in \cite{DOR03}, the reduction of contact Lagrangian systems (Lagrange--Poincar\'e--Herglotz equations) in \cite{ACdLSS24}, contact mechanical systems on Lie groups in \cite{ACdLMMdDP}, momentum maps and contact reduction, together with their relation with symplectization, in \cite{MGM22}, and reduction by scaling symmetries leading to contact systems in \cite{BJS23}; coisotropic reduction in contact geometry was established in \cite{dLL19} and is reviewed in \cite{dLIL24}. The aim of this paper is to develop, for contact Hamiltonian systems, the analogue of the theory of \cite{dLMdDV17}, that is, to understand how the four contact Hamilton--Jacobi equations interact with symmetries, momentum maps and coisotropic reduction.

\medskip\noindent\textbf{Main results.} Let us summarize our contributions.

\begin{enumerate}[leftmargin=*,label=(\arabic*)]
\item \emph{Geometry of the contact Hamilton--Jacobi equations} (Section \ref{sec:HJ}). We show that the image of a section $\gamma\colon Q\times\R\to T^*Q\times\R$ is coisotropic if and only if the $1$-form $\gamma^*\eta_Q$ is integrable, if and only if $\Im\gamma$ is foliated by Legendrian submanifolds; locally these leaves are the $1$-jets $j^1S_c$ of a one-parameter family of functions (Proposition \ref{prop:coiso-sections}). The additional assumption of Lagrangian slices used in \cite{dLLMB21} is not needed. With this description, the Hamilton--Jacobi equation for $\mathcal{E}_H$ means that $H$ is constant on each Legendrian leaf (Theorem \ref{thm:case3}), and the one for $X_H$ is equivalent to the evolutionary contact Hamilton--Jacobi equation $\partial_tS+H(q,\partial_qS,S)=0$, the leaves being transported by the flow of $X_H$ (Theorem \ref{thm:case1} and Corollary \ref{cor:evolution}). Legendrian solutions are characterized in a coordinate-free way by means of the elementary Lemma \ref{lem:tangency} (compare with \cite{dLL19,dLLMB21} and \cite{ELdLM21}): a Legendrian submanifold is invariant under $X_f$ (resp.\ $\mathcal{E}_f$) if and only if $f$ vanishes (resp.\ is locally constant) on it. This allows us to define generalized (non-horizontal) solutions.

\item \emph{Momentum maps and invariant Legendrian submanifolds} (Section \ref{sec:symmetries}). For an action by contactomorphisms with momentum map $J$, we characterize the regular values of $J$ (Proposition \ref{prop:regular}), we prove that $J$ is dissipated along the flows of $X_H$ and $\mathcal{E}_H$, so that $J^{-1}(0)$ is the only level set that is invariant under the dynamics of every $G$-invariant Hamiltonian $H$ with $\Reeb(H)\neq0$ (Proposition \ref{prop:dissipated}); in particular, if $\Reeb(H)$ does not vanish, there are no Legendrian solutions of the Hamilton--Jacobi problem on which $J$ takes a constant non-zero value (Corollary \ref{cor:no-nonzero}), and we prove the contact analogue of the key lemma of \cite{dLMdDV17}: a Legendrian submanifold is $G$-invariant if and only if $J$ vanishes on it; if $J$ is constant and equal to $\mu$, then $\mu$ vanishes on $[\g,\g]$ (equivalently, $\Ad^*_g\mu=\mu$ for all $g\in G$) and the submanifold is invariant under the twisted fields $\xi_M+\mu(\xi)\Reeb$, and conversely (Theorem \ref{thm:invariant-legendrian}). The Hamilton--Jacobi equations for Legendrian sections appear as a particular case. Finally, for a regular value $\mu$, $J^{-1}(\mu)$ is coisotropic if and only if $\mu$ vanishes on $[\g,\g]$, and its characteristic distribution is generated by the twisted fields (Proposition \ref{prop:levelsets}).

\item \emph{Contact coisotropic reduction} (Section \ref{sec:reduction}). We recall, with self-contained proofs adapted to our needs, the contact coisotropic reduction theorem and the reduction of Legendrian submanifolds (see \cite{dLL19} and \cite[Theorem 7.1 and Proposition 7.7]{dLIL24}): if $C$ is a coisotropic submanifold tangent to the Reeb vector field, its characteristic distribution is the kernel of the pair $(i_C^*\eta,i_C^*d\eta)$ and the leaf space inherits a contact form; Legendrian submanifolds with clean intersection with $C$ project onto Legendrian submanifolds (Theorems \ref{thm:coiso-red} and \ref{thm:legendrian-red}). For cotangent lifted actions, the contact reduced space at zero is $(T^*(Q/G)\times\R,\eta_{Q/G})$, and the reduced dynamics is the contact dynamics of the reduced Hamiltonian (Theorem \ref{thm:zero-reduction} and Proposition \ref{prop:cotangent-red}).

\item \emph{Reduction and reconstruction of the Hamilton--Jacobi equation} (Section \ref{sec:HJ-reduction}). We establish a one-to-one correspondence between $G$-invariant (generalized) solutions of the contact Hamilton--Jacobi problem for $H$ and solutions for the reduced Hamiltonian, for the four equations (Theorems \ref{thm:HJ-red}, \ref{thm:HJ-rec} and \ref{thm:HJ-red-coiso}), and we describe the reconstruction of the dynamics (Proposition \ref{prop:reconstruction}). In contrast with \cite{dLMdDV17}, $G$-invariant solutions only occur at the zero level set of the momentum map, no magnetic terms appear, and invariant Legendrian solutions always have invariant generating functions.
\item \emph{Reduction at non-zero momentum} (Section \ref{sec:nonzero}). For $\mu\neq0$ vanishing on $[\g,\g]$ and integrating to a character $\chi\colon G\to\R$, we construct a strict contactomorphism $T_\varphi$ of $T^*Q\times\R$ that maps $J^{-1}(\mu)$ onto $J^{-1}(0)$ and conjugates the twisted action to the original one (Proposition \ref{prop:shift}). This yields a reduction of the Hamilton--Jacobi problem at $\mu$ when $\Reeb(H)=0$ on $J^{-1}(\mu)$, which contains the classical ansatz for cyclic coordinates (Theorem \ref{thm:HJ-twisted}). For arbitrary dissipative Hamiltonians, the preimage $J^{-1}(\R^+\mu)$ of the ray is invariant; it is reduced by the subgroup $K$ with Lie algebra $\ker\mu$ to an open subset of $T^*(Q/K)\times\R$, which is Willett's reduced space, and solutions with momentum on the ray correspond to solutions of a reduced equation with $\dim G-1$ fewer variables; this is the reduction at zero by $K$, restricted to an open set (Proposition \ref{prop:ray} and Theorem \ref{thm:HJ-ray}).
\end{enumerate}

Section \ref{sec:examples} contains examples: cyclic coordinates in the presence of dissipation, a central force problem with linear friction, a two-body problem with linear dissipation and a particle in the plane with friction, which illustrates the reduction along rays; in the dissipative case the classical ansatz for cyclic coordinates only survives at zero momentum. We finish with some conclusions and open problems in Section \ref{sec:conclusions}.

Throughout the paper all manifolds and maps are smooth, manifolds are connected unless otherwise stated, and the flows of the vector fields are assumed to be complete (all the results hold locally with the obvious modifications). Lie groups are assumed to be connected. A submanifold is called \emph{closed} if it is embedded and closed as a subset.

\section{Preliminaries on contact geometry}\label{sec:prelim}

\subsection{Contact manifolds and contact Hamiltonian systems}

\begin{definition}
A \emph{contact manifold} is a pair $(M,\eta)$ where $M$ is a $(2n+1)$-dimensional manifold and $\eta$ is a $1$-form on $M$ such that $\eta\wedge(d\eta)^n$ is a volume form. The $1$-form $\eta$ is called a \emph{contact form}.
\end{definition}

Given a contact manifold $(M,\eta)$ there is a unique vector field $\Reeb$ on $M$, the \emph{Reeb vector field}, such that
\[
i_\Reeb d\eta=0,\qquad i_\Reeb\eta=1.
\]
We have the Whitney sum decomposition $TM=\Hor\oplus\Ver$, where $\Hor=\ker\eta$ is the \emph{horizontal} (or contact) distribution, of rank $2n$, and $\Ver=\ker d\eta=\langle\Reeb\rangle$ is the \emph{vertical} distribution, of rank one. For each $x\in M$, $(\Hor_x,d\eta_x)$ is a symplectic vector space. If $W\subseteq\Hor_x$ is a linear subspace we denote by
\[
W^{\perpd}=\{v\in\Hor_x\mid d\eta_x(v,w)=0\ \text{for all } w\in W\}
\]
its symplectic orthogonal \emph{inside} $\Hor_x$.

By the Darboux theorem \cite{Gei08,LM87}, around any point there are local coordinates $(q^i,p_i,z)$, $i=1,\dots,n$, such that
\[
\eta=dz-p_i\,dq^i,\qquad d\eta=dq^i\wedge dp_i,\qquad \Reeb=\pd{}{z}.
\]
In these coordinates, $\Hor$ is spanned by the vector fields
\begin{equation}\label{eq:AB}
A_i=\pd{}{q^i}+p_i\pd{}{z},\qquad B^i=\pd{}{p_i},\qquad d\eta(A_i,B^j)=\delta_i^j .
\end{equation}

The contact form defines a vector bundle isomorphism
\[
\flat\colon TM\to T^*M,\qquad \flat(v)=i_vd\eta+\eta(v)\,\eta,
\]
whose inverse we denote by $\sharp$. Notice that $\flat(\Reeb)=\eta$ and $\flat(v)(\Reeb)=\eta(v)$. We use the same notation for the induced isomorphism of $C^\infty(M)$-modules $\X(M)\to\Omega^1(M)$.

\begin{definition}
A \emph{diffeomorphism} $F\colon(M,\eta)\to(N,\xi)$ is a \emph{contactomorphism} if $F^*\xi=\eta$, and a \emph{conformal contactomorphism} if $F^*\xi=a\,\eta$ for some nowhere vanishing function $a$. A vector field $X$ is an \emph{infinitesimal contactomorphism} (resp.\ \emph{infinitesimal conformal contactomorphism}) if $\Lie_X\eta=0$ (resp.\ $\Lie_X\eta=b\,\eta$ for some function $b$).
\end{definition}

A contactomorphism $F$ maps the Reeb vector field of $\eta$ onto the Reeb vector field of $\xi$, by uniqueness. In particular, if $X$ is an infinitesimal contactomorphism then $[X,\Reeb]=0$.

\begin{definition}
Let $H\in C^\infty(M)$. The \emph{Hamiltonian vector field} of $H$ is the vector field $X_H$ defined by
\[
\flat(X_H)=dH-\big(\Reeb(H)+H\big)\eta,
\]
or equivalently by the two conditions
\begin{equation}\label{eq:XH}
\eta(X_H)=-H,\qquad i_{X_H}d\eta=dH-\Reeb(H)\,\eta .
\end{equation}
The \emph{evolution vector field} of $H$ is $\mathcal{E}_H=X_H+H\Reeb$; it is characterized by
\begin{equation}\label{eq:EH}
\eta(\mathcal{E}_H)=0,\qquad i_{\mathcal{E}_H}d\eta=dH-\Reeb(H)\,\eta .
\end{equation}
The triple $(M,\eta,H)$ is called a \emph{contact Hamiltonian system}.
\end{definition}

In Darboux coordinates,
\begin{align}
X_H&=\pd{H}{p_i}\pd{}{q^i}-\left(\pd{H}{q^i}+p_i\pd{H}{z}\right)\pd{}{p_i}+\left(p_i\pd{H}{p_i}-H\right)\pd{}{z},\label{eq:XHloc}\\
\mathcal{E}_H&=\pd{H}{p_i}\pd{}{q^i}-\left(\pd{H}{q^i}+p_i\pd{H}{z}\right)\pd{}{p_i}+p_i\pd{H}{p_i}\pd{}{z}.\label{eq:EHloc}
\end{align}
Hence the integral curves of $X_H$ are the solutions of the \emph{contact Hamilton equations}
\[
\dot q^i=\pd{H}{p_i},\qquad \dot p_i=-\left(\pd{H}{q^i}+p_i\pd{H}{z}\right),\qquad \dot z=p_i\pd{H}{p_i}-H,
\]
which are the Herglotz equations for the Legendre transform of a regular Lagrangian \cite{He30,dLL19}, while the integral curves of $\mathcal{E}_H$ satisfy the same equations with $\dot z=p_i\,\partial H/\partial p_i$.

\begin{proposition}\label{prop:basic-dyn}
For every $H\in C^\infty(M)$ we have
\[
\Lie_{X_H}\eta=-\Reeb(H)\,\eta,\qquad X_H(H)=-\Reeb(H)\,H,\qquad \mathcal{E}_H(H)=0 .
\]
In particular, $X_H$ is an infinitesimal conformal contactomorphism, which is an infinitesimal contactomorphism if and only if $\Reeb(H)=0$.
\end{proposition}

\begin{proof}
By \eqref{eq:XH}, $\Lie_{X_H}\eta=d(\eta(X_H))+i_{X_H}d\eta=-dH+dH-\Reeb(H)\eta$. Moreover $X_H(H)=dH(X_H)=d\eta(X_H,X_H)+\Reeb(H)\eta(X_H)=-\Reeb(H)H$, and $\mathcal{E}_H(H)=d\eta(\mathcal{E}_H,\mathcal{E}_H)+\Reeb(H)\eta(\mathcal{E}_H)=0$ by \eqref{eq:EH}.
\end{proof}

Conversely, every infinitesimal conformal contactomorphism $X$ is the Hamiltonian vector field of the function $f=-\eta(X)$; see \cite{dLL19,dLL20}. The evolution vector field plays a relevant role in thermodynamics \cite{SdLLMdD20,dLB25}, since $H$ is a conserved quantity for it.

\subsection{Contact manifolds as Jacobi manifolds}

\begin{definition}[\cite{Ki76,Li78}]
A \emph{Jacobi manifold} is a triple $(M,\Lambda,Z)$ where $\Lambda$ is a bivector field and $Z$ is a vector field such that
\[
[\Lambda,\Lambda]=2Z\wedge\Lambda,\qquad \Lie_Z\Lambda=[Z,\Lambda]=0,
\]
$[\cdot,\cdot]$ being the Schouten--Nijenhuis bracket. The \emph{Jacobi bracket} is $\{f,g\}=\Lambda(df,dg)+fZ(g)-gZ(f)$.
\end{definition}

The Jacobi bracket is bilinear, skew-symmetric and satisfies the Jacobi identity, but instead of the Leibniz rule it satisfies the weak Leibniz rule $\{f,gh\}=g\{f,h\}+h\{f,g\}+ghZ(f)$. A Poisson manifold is a Jacobi manifold with $Z=0$.

A contact manifold $(M,\eta)$ is a Jacobi manifold with
\[
\Lambda(\alpha,\beta)=-d\eta(\sharp\alpha,\sharp\beta),\qquad Z=-\Reeb .
\]
Let $\sharpL\colon T^*M\to TM$ be defined by $\langle\sharpL(\alpha),\beta\rangle=\Lambda(\alpha,\beta)$. In Darboux coordinates $\Lambda=\partial_{p_i}\wedge\big(\partial_{q^i}+p_i\partial_z\big)$ and $\sharpL(\alpha)=\alpha(B^i)A_i-\alpha(A_i)B^i$. One easily checks the intrinsic expression
\begin{equation}\label{eq:sharpL}
\sharpL(\alpha)=\sharp(\alpha)-\alpha(\Reeb)\,\Reeb .
\end{equation}

\begin{lemma}\label{lem:sharpL}
The map $\sharpL$ has kernel $\langle\eta\rangle$ and image $\Hor$. Moreover:
\begin{enumerate}[label=(\roman*)]
\item $d\eta(\sharpL(\alpha),v)=\alpha(v)$ for every $\alpha\in T_x^*M$ and $v\in\Hor_x$;
\item $\sharpL(df)=\mathcal{E}_f$ and hence $X_f=\sharpL(df)-f\Reeb$, for every $f\in C^\infty(M)$.
\end{enumerate}
\end{lemma}

\begin{proof}
Since $\flat(v)(\Reeb)=\eta(v)$, we have $\eta(\sharp\alpha)=\alpha(\Reeb)$, so $\eta(\sharpL\alpha)=0$ by \eqref{eq:sharpL}. If $\sharpL\alpha=0$, then $\sharp\alpha=\alpha(\Reeb)\Reeb$ and $\alpha=\alpha(\Reeb)\eta$; by dimension, the image is $\Hor$. For $v\in\Hor_x$, $d\eta(\sharpL\alpha,v)=d\eta(\sharp\alpha,v)=\flat(\sharp\alpha)(v)-\eta(\sharp\alpha)\eta(v)=\alpha(v)$. Finally, $\sharp(dH)=X_H+(\Reeb(H)+H)\Reeb$ by definition of $X_H$, so $\sharpL(dH)=X_H+H\Reeb=\mathcal{E}_H$.
\end{proof}

\subsection{Submanifolds of a contact manifold}

\begin{definition}[\cite{dLL19,dLLMB21}]
Let $\Delta_x\subseteq T_xM$ be a linear subspace. Its \emph{contact complement} is $\Delta_x^{\perpL}=\sharpL(\Delta_x^\circ)$, where $\Delta_x^\circ\subseteq T_x^*M$ is the annihilator of $\Delta_x$. A submanifold $N\subseteq M$ is said to be
\begin{itemize}
\item \emph{isotropic} if $TN\subseteq TN^{\perpL}$;
\item \emph{coisotropic} if $TN^{\perpL}\subseteq TN$;
\item \emph{Legendrian} if $TN=TN^{\perpL}$.
\end{itemize}
\end{definition}

The following lemma reduces the contact complement to the symplectic orthogonal in $(\Hor,d\eta)$ and is the main tool in what follows.

\begin{lemma}\label{lem:complement}
For every linear subspace $\Delta_x\subseteq T_xM$ we have
\[
\Delta_x^{\perpL}=(\Delta_x\cap\Hor_x)^{\perpd}.
\]
\end{lemma}

\begin{proof}
By Lemma \ref{lem:sharpL}, $\Delta_x^{\perpL}\subseteq\Hor_x$, and for $\alpha\in\Delta_x^\circ$ and $w\in\Delta_x\cap\Hor_x$ we get $d\eta(\sharpL\alpha,w)=\alpha(w)=0$. Hence $\Delta_x^{\perpL}\subseteq(\Delta_x\cap\Hor_x)^{\perpd}$. Let $m=\dim\Delta_x$. If $\Delta_x\not\subseteq\Hor_x$, then $\eta_x\notin\Delta_x^\circ$, $\sharpL$ is injective on $\Delta_x^\circ$ and both spaces have dimension $2n+1-m$. If $\Delta_x\subseteq\Hor_x$, then $\eta_x\in\Delta_x^\circ$ and both spaces have dimension $2n-m$.
\end{proof}

As a consequence: a submanifold $N$ is isotropic if and only if $TN\subseteq\Hor$ (that is, $N$ is an integral submanifold of $\Hor$; note that then $d\eta|_{TN}=0$ automatically); $N$ is coisotropic if and only if $(T_xN\cap\Hor_x)^{\perpd}\subseteq T_xN$ for all $x\in N$; and we recover the following well-known characterization.

\begin{proposition}[{\cite{dLL19,Gei08}}]\label{prop:legendrian}
A submanifold $L$ of a $(2n+1)$-dimensional contact manifold is Legendrian if and only if it is an integral submanifold of $\Hor$ of maximal dimension $n$, that is, $\dim L=n$ and $i_L^*\eta=0$. Equivalently, $T_xL$ is a Lagrangian subspace of $(\Hor_x,d\eta_x)$ for every $x\in L$.
\end{proposition}

In local coordinates, if $N$ is the zero set of independent functions $\phi_a$, then $TN^{\perpL}$ is spanned by $\sharpL(d\phi_a)=\mathcal{E}_{\phi_a}$, and $N$ is coisotropic if and only if
\begin{equation}\label{eq:coiso-local}
\Lambda(d\phi_a,d\phi_b)=\mathcal{E}_{\phi_a}(\phi_b)=\pd{\phi_a}{p_i}\left(\pd{\phi_b}{q^i}+p_i\pd{\phi_b}{z}\right)-\left(\pd{\phi_a}{q^i}+p_i\pd{\phi_a}{z}\right)\pd{\phi_b}{p_i}=0\quad\text{on } N.
\end{equation}

The following characterization of coisotropic submanifolds of dimension $n+1$ will be essential; it is closely related to the results on coisotropic and Legendre--Lagrangian submanifolds of Jacobi manifolds in \cite{IdLMdD97}.

\begin{lemma}\label{lem:coiso-foliation}
Let $N\subseteq M$ be a submanifold of dimension $n+1$ which is transverse to $\Hor$ (i.e., $T_xN\not\subseteq\Hor_x$ for all $x\in N$; this always holds on an open dense subset of $N$, since integral submanifolds of $\Hor$ have dimension at most $n$), and let $\vartheta_N=i_N^*\eta$, a nowhere vanishing $1$-form on $N$. The following assertions are equivalent:
\begin{enumerate}[label=(\roman*)]
\item $N$ is coisotropic;
\item $\vartheta_N\wedge d\vartheta_N=0$, that is, the distribution $\ker\vartheta_N=TN\cap\Hor$ is involutive;
\item $N$ is foliated by Legendrian submanifolds of $M$.
\end{enumerate}
In this case, the Legendrian leaves are the leaves of the foliation $\ker\vartheta_N$.
\end{lemma}

\begin{proof}
Since $N$ is transverse to $\Hor$, $T_xN\cap\Hor_x$ has dimension $n$, so $(T_xN\cap\Hor_x)^{\perpd}$ has dimension $n$ too, and it is contained in $\Hor_x$. Hence, by Lemma \ref{lem:complement}, $N$ is coisotropic if and only if $(T_xN\cap\Hor_x)^{\perpd}=T_xN\cap\Hor_x$, i.e., if and only if $T_xN\cap\Hor_x$ is a Lagrangian subspace of $(\Hor_x,d\eta_x)$, which means that $d\vartheta_N$ vanishes on $\ker\vartheta_N$. This is equivalent to $\vartheta_N\wedge d\vartheta_N=0$ and, by the Frobenius theorem, to the involutivity of $\ker\vartheta_N$. If this holds, the leaves of $\ker\vartheta_N$ are $n$-dimensional integral submanifolds of $\Hor$, that is, Legendrian submanifolds. Conversely, if $N$ is foliated by Legendrian submanifolds, the tangent spaces of the leaves are $n$-dimensional subspaces of $TN\cap\Hor$, hence they coincide with it, and $\ker\vartheta_N$ is involutive.
\end{proof}

\subsection{The canonical contact manifold $T^*Q\times\R$ and $1$-jets}

Let $Q$ be an $n$-dimensional manifold, $\theta_Q$ the Liouville $1$-form on $T^*Q$ and $\omega_Q=-d\theta_Q$ the canonical symplectic form. The manifold $T^*Q\times\R$ is equipped with the canonical contact form
\[
\eta_Q=dz-\rho^*\theta_Q,
\]
where $\rho\colon T^*Q\times\R\to T^*Q$ is the projection and $z$ is the coordinate in $\R$. Bundle coordinates $(q^i,p_i,z)$ are Darboux coordinates, and $d\eta_Q=\rho^*\omega_Q$. We will use the projections
\[
\pi\colon T^*Q\times\R\to Q\times\R,\qquad \pi_1\colon T^*Q\times\R\to Q,\qquad \pi_Q\colon T^*Q\to Q .
\]
The $1$-jet of a function $f\colon Q\to\R$ is the section $j^1f\colon Q\to T^*Q\times\R$ of $\pi_1$ given by $j^1f(q)=(df(q),f(q))$.

\begin{proposition}\label{prop:jets}
A section $\gamma=(\gamma_Q,\gamma_z)\colon Q\to T^*Q\times\R$ of $\pi_1$ has Legendrian image if and only if $\gamma=j^1f$ with $f=\gamma_z$. In particular, $\gamma(Q)$ is Legendrian if and only if $\gamma$ is the $1$-jet of a function, and then $\rho\circ\gamma=df$ has Lagrangian image in $(T^*Q,\omega_Q)$.
\end{proposition}

\begin{proof}
Since $\gamma^*\rho^*\theta_Q=\gamma_Q$ (the Liouville form is tautological), $\gamma^*\eta_Q=d\gamma_z-\gamma_Q$, which vanishes if and only if $\gamma_Q=d\gamma_z$. The dimension of $\gamma(Q)$ is $n$.
\end{proof}

This is the contact analogue of the well-known fact that the image of a $1$-form on $Q$ is Lagrangian in $(T^*Q,\omega_Q)$ if and only if the form is closed.

\subsection{Legendrian submanifolds and contact vector fields}

The following simple lemma, which is essentially known (see \cite{dLL19,dLLMB21} and \cite{ELdLM21}), will be used repeatedly; we include a short proof for completeness. It gives an intrinsic version of the Hamilton--Jacobi equations for Legendrian submanifolds and will be the key to relate them with momentum maps.

\begin{lemma}\label{lem:tangency}
Let $L\subseteq M$ be a connected Legendrian submanifold and $f\in C^\infty(M)$.
\begin{enumerate}[label=(\roman*)]
\item $X_f$ is tangent to $L$ if and only if $f|_L=0$.
\item $\mathcal{E}_f$ is tangent to $L$ if and only if $f|_L$ is constant. In this case, $X_f+(f|_L)\,\Reeb$ is tangent to $L$.
\end{enumerate}
\end{lemma}

\begin{proof}
For $v\in T_xL\subseteq\Hor_x$, equations \eqref{eq:XH} and \eqref{eq:EH} give
\begin{equation}\label{eq:aux-tangency}
d\eta(X_f,v)=d\eta(\mathcal{E}_f,v)=df(v).
\end{equation}
(i) If $X_f$ is tangent to $L$, then $f=-\eta(X_f)=0$ on $L$. Conversely, if $f|_L=0$ then $X_f(x)\in\Hor_x$ and, by \eqref{eq:aux-tangency}, $X_f(x)\in(T_xL)^{\perpd}=T_xL$, because $T_xL$ is Lagrangian in $\Hor_x$.
(ii) Since $\eta(\mathcal{E}_f)=0$, $\mathcal{E}_f(x)\in\Hor_x$. By \eqref{eq:aux-tangency}, $\mathcal{E}_f(x)\in(T_xL)^{\perpd}=T_xL$ if and only if $df$ vanishes on $T_xL$. The last assertion follows from $\mathcal{E}_f=X_f+f\Reeb$.
\end{proof}

\section{The Hamilton--Jacobi equations for contact Hamiltonian systems}\label{sec:HJ}

In this section we revisit the four Hamilton--Jacobi equations obtained in \cite{dLLMB21} and give a geometric interpretation of their solutions which will be crucial in the study of symmetries. Throughout this section $H\colon T^*Q\times\R\to\R$ is a Hamiltonian function.

The following elementary lemma explains the meaning of the Hamilton--Jacobi condition.

\begin{lemma}\label{lem:related-tangent}
Let $\tau\colon P\to B$ be a fibre bundle, $\gamma\colon B\to P$ a section and $X\in\X(P)$. Define the projected vector field $X^\gamma=T\tau\circ X\circ\gamma\in\X(B)$. Then $X\circ\gamma=T\gamma\circ X^\gamma$ (that is, $X$ and $X^\gamma$ are $\gamma$-related) if and only if $X$ is tangent to the submanifold $\gamma(B)$.
\end{lemma}

\begin{proof}
Since $\gamma$ is an embedding, $T_{\gamma(b)}\gamma(B)=\Im T_b\gamma$. If $X(\gamma(b))=T_b\gamma(v)$, applying $T\tau$ and using $\tau\circ\gamma=\id_B$ we get $v=X^\gamma(b)$. The converse is obvious.
\end{proof}

\subsection{Coisotropic sections and Legendrian foliations}

Let $\gamma\colon Q\times\R\to T^*Q\times\R$ be a section of $\pi$, locally $\gamma(q,z)=(q^i,\gamma_i(q,z),z)$. We denote by $\bar\gamma=\gamma_i\,dq^i$ the corresponding $z$-dependent $1$-form on $Q$ and introduce
\[
\vartheta_\gamma:=\gamma^*\eta_Q=dz-\gamma_i\,dq^i\in\Omega^1(Q\times\R),\qquad A_j^\gamma=\pd{}{q^j}+\gamma_j\pd{}{z}\in\X(Q\times\R),
\]
so that $\ker\vartheta_\gamma=\langle A_1^\gamma,\dots,A_n^\gamma\rangle$ and $\vartheta_\gamma(\partial_z)=1$. The submanifold $N_\gamma=\gamma(Q\times\R)$ has dimension $n+1$ and it is transverse to $\Hor$, since $\eta_Q(T\gamma(\partial_z))=1$.

\begin{proposition}\label{prop:coiso-sections}
The following assertions are equivalent:
\begin{enumerate}[label=(\roman*)]
\item $N_\gamma=\gamma(Q\times\R)$ is a coisotropic submanifold of $(T^*Q\times\R,\eta_Q)$;
\item $\vartheta_\gamma\wedge d\vartheta_\gamma=0$;
\item $[A_i^\gamma,A_j^\gamma]=0$ for all $i,j$, that is,
\begin{equation}\label{eq:coiso-gamma}
A_i^\gamma(\gamma_j)=A_j^\gamma(\gamma_i),\quad\text{i.e.}\quad \pd{\gamma_j}{q^i}+\gamma_i\pd{\gamma_j}{z}=\pd{\gamma_i}{q^j}+\gamma_j\pd{\gamma_i}{z};
\end{equation}
\item $N_\gamma$ is foliated by Legendrian submanifolds;
\item locally, there is a function $S(q,c)$ with $\partial S/\partial c\neq0$ such that
\begin{equation}\label{eq:family}
\gamma_j\big(q,S(q,c)\big)=\pd{S}{q^j}(q,c);
\end{equation}
that is, $N_\gamma=\bigcup_c j^1S_c(Q)$ with $S_c=S(\cdot,c)$, and the Legendrian leaves of $N_\gamma$ are the $1$-jets $j^1S_c(Q)$.
\end{enumerate}
\end{proposition}

\begin{proof}
The equivalence of (i), (ii) and (iv) is Lemma \ref{lem:coiso-foliation}, since $\gamma\colon Q\times\R\to N_\gamma$ is a diffeomorphism and $\gamma^*(i_{N_\gamma}^*\eta_Q)=\vartheta_\gamma$. Since $[A_i^\gamma,A_j^\gamma]=\big(A_i^\gamma(\gamma_j)-A_j^\gamma(\gamma_i)\big)\partial_z$ is transverse to $\ker\vartheta_\gamma$ unless it vanishes, (ii) and (iii) are equivalent by the Frobenius theorem. If (ii) holds, locally $\vartheta_\gamma=\lambda\,du$ for some functions $\lambda,u$ with $\lambda\,\partial u/\partial z=1$. By the implicit function theorem, $u(q,z)=c$ if and only if $z=S(q,c)$ with $\partial S/\partial c=1/(\partial u/\partial z)\neq0$, and since $\vartheta_\gamma$ vanishes along the leaves $\{u=c\}$ we get $dS_c=\gamma_j(q,S_c(q))\,dq^j$, which is \eqref{eq:family}. Conversely, given $S$ as in (v), in the local coordinates $(q,c)$ of $Q\times\R$ defined by $z=S(q,c)$ we have $\vartheta_\gamma=dS-\partial_jS\,dq^j=\pd{S}{c}\,dc$, which is integrable.
\end{proof}

\begin{remark}\label{rem:lagrangian-slices}
In \cite{dLLMB21} the Hamilton--Jacobi problem for sections of $\pi$ is formulated under two assumptions: $N_\gamma$ is coisotropic and the slices $\gamma^z=\rho\circ\gamma(\cdot,z)$ have Lagrangian image in $(T^*Q,\omega_Q)$, i.e., $d_Q\bar\gamma=0$, where $d_Q$ denotes the exterior differential at fixed $z$. Under \eqref{eq:coiso-gamma}, the last condition is equivalent to $\bar\gamma\wedge\partial_z\bar\gamma=0$, that is, $\partial_z\bar\gamma=\sigma\bar\gamma$ for some function $\sigma$ on the open set where $\bar\gamma\neq0$ (compare with \cite[Theorem 2]{dLLMB21}). As we will see, \emph{only coisotropy is needed} to obtain the Hamilton--Jacobi equations, and the natural objects are the Legendrian leaves $j^1S_c(Q)$, which are in general different from the slices $\gamma(Q\times\{z\})$.
\end{remark}

\subsection{Case 1: the Hamiltonian vector field and coisotropic sections}

Let $\gamma$ be a section of $\pi$ and $X_H^\gamma=T\pi\circ X_H\circ\gamma$. Write $F_\gamma=H\circ\gamma\in C^\infty(Q\times\R)$.

\begin{theorem}\label{thm:case1}
Assume that $N_\gamma=\gamma(Q\times\R)$ is coisotropic. Then $X_H$ and $X_H^\gamma$ are $\gamma$-related if and only if
\begin{equation}\label{eq:HJ1}
A_j^\gamma(F_\gamma)=F_\gamma\,\pd{\gamma_j}{z},\qquad j=1,\dots,n,
\end{equation}
or, intrinsically,
\begin{equation}\label{eq:HJ1-intr}
dF_\gamma-\pd{F_\gamma}{z}\,\vartheta_\gamma+F_\gamma\, i_{\partial_z}d\vartheta_\gamma=0 .
\end{equation}
Equivalently, \eqref{eq:HJ1-intr} can be written as
\[
d_Q(H\circ\gamma)+\gamma_o\,\bar\gamma-(H\circ\gamma)\,i_{\partial_z}d(\gamma^*\theta_Q)=0,\qquad \gamma_o=\pd{H}{z}+\pd{H}{p_i}\pd{\gamma_i}{z}=\pd{(H\circ\gamma)}{z},
\]
which is the Hamilton--Jacobi equation of \cite[Theorem 1]{dLLMB21}, obtained there under the additional assumption that the slices $\gamma^z$ are Lagrangian; the identity itself does not require this assumption.
\end{theorem}

\begin{proof}
From \eqref{eq:XHloc}, $X_H^\gamma=\pd{H}{p_i}\partial_{q^i}+\big(\gamma_i\pd{H}{p_i}-H\big)\partial_z$, the partial derivatives of $H$ being evaluated at $\gamma(q,z)$. Then
\[
T\gamma(X_H^\gamma)=\pd{H}{p_i}\pd{}{q^i}+\left(\gamma_i\pd{H}{p_i}-H\right)\left(\pd{}{z}+\pd{\gamma_j}{z}\pd{}{p_j}\right)+\pd{H}{p_i}\pd{\gamma_j}{q^i}\pd{}{p_j}.
\]
The $\partial_{q}$ and $\partial_z$ components agree with those of $X_H\circ\gamma$, so the $\gamma$-relatedness is equivalent to the equality of the $\partial_{p_j}$ components:
\begin{equation}\label{eq:case1-coord}
\pd{H}{q^j}+\gamma_j\pd{H}{z}+\pd{H}{p_i}A_i^\gamma(\gamma_j)-H\pd{\gamma_j}{z}=0 .
\end{equation}
By the chain rule, $A_j^\gamma(F_\gamma)=\pd{H}{q^j}+\gamma_j\pd{H}{z}+\pd{H}{p_i}A_j^\gamma(\gamma_i)$. Using the coisotropy condition \eqref{eq:coiso-gamma}, \eqref{eq:case1-coord} becomes \eqref{eq:HJ1}. Finally, $d\vartheta_\gamma=-d\gamma_i\wedge dq^i$, so $i_{\partial_z}d\vartheta_\gamma=-\pd{\gamma_i}{z}dq^i$ and the left-hand side of \eqref{eq:HJ1-intr} is $\big(A^\gamma_j(F_\gamma)-F_\gamma\pd{\gamma_j}{z}\big)dq^j$. The last assertion follows from $\gamma^*\theta_Q=\bar\gamma$ and $d_Q(H\circ\gamma)+\pd{(H\circ\gamma)}{z}\bar\gamma=A_j^\gamma(F_\gamma)dq^j$.
\end{proof}

A section $\gamma$ with coisotropic image satisfying \eqref{eq:HJ1} is called a \emph{solution of the Hamilton--Jacobi problem for $X_H$}. The meaning of \eqref{eq:HJ1} is clarified in terms of the Legendrian leaves.

\begin{proposition}\label{prop:case1-leaves}
Let $\gamma$ be a section with coisotropic image, locally described by a function $S(q,c)$ as in Proposition \ref{prop:coiso-sections}, and let $h(q,c)=H\big(j^1S_c(q)\big)=H\big(q,\partial_qS(q,c),S(q,c)\big)$. Then $\gamma$ is a solution of the Hamilton--Jacobi problem for $X_H$ if and only if there is a function $k$ of $c$ alone such that
\begin{equation}\label{eq:HJ1-leaves}
H\left(q,\pd{S}{q}(q,c),S(q,c)\right)=k(c)\,\pd{S}{c}(q,c).
\end{equation}
\end{proposition}

\begin{proof}
Since $F_\gamma(q,S(q,c))=h(q,c)$ and, by \eqref{eq:family}, $\partial_{q^j}\big(S(q,c)\big)=\gamma_j(q,S(q,c))$, the chain rule gives $\partial h/\partial q^j=A_j^\gamma(F_\gamma)$ at $(q,S(q,c))$. Differentiating \eqref{eq:family} with respect to $c$ gives $\partial_z\gamma_j\,\partial_cS=\partial_j\partial_cS$. Hence \eqref{eq:HJ1} reads $\partial_jh\,\partial_cS-h\,\partial_j\partial_cS=0$, that is, $\partial_j\big(h/\partial_cS\big)=0$ for all $j$.
\end{proof}

\begin{corollary}\label{cor:evolution}
Let $\gamma$ be a solution of the Hamilton--Jacobi problem for $X_H$ and assume that $H\circ\gamma\neq0$. Then the Legendrian leaves of $N_\gamma$ can be (locally) parametrized by a parameter $t$ such that the corresponding function $S(q,t)$ satisfies the \emph{evolutionary contact Hamilton--Jacobi equation}
\begin{equation}\label{eq:evolHJ}
\pd{S}{t}+H\left(q,\pd{S}{q},S\right)=0 .
\end{equation}
Moreover, if $\varphi_s$ denotes the flow of $X_H$, then $\varphi_s\big(j^1S_t(Q)\big)=j^1S_{t+s}(Q)$. Conversely, every solution $S(q,t)$ of \eqref{eq:evolHJ} with $H(q,\partial_qS,S)\neq0$ defines, by \eqref{eq:family}, a solution of the Hamilton--Jacobi problem for $X_H$.
\end{corollary}

\begin{proof}
Since $h=k(c)\partial_cS$ and $h\neq0$, $k$ does not vanish. Define $t$ by $dc/dt=-k(c)$; then $\partial_tS=\partial_cS\,dc/dt=-h$, which is \eqref{eq:evolHJ}. Consider the parametrization $\Psi(q,t)=j^1S_t(q)$ of $N_\gamma$. Since $X_H$ is tangent to $N_\gamma$, $X_H(\Psi)=T\Psi(a\,\partial_t+b^j\partial_{q^j})$; applying $\eta_Q$ and using that $\eta_Q$ vanishes on the leaves, $-H=a\,\partial_tS=-a\,H$, so $a=1$, which proves the assertion on the flow. For the converse, note that $\partial_tS=-H(q,\partial_qS,S)\neq0$, so that $S$ satisfies the non-degeneracy condition of Proposition \ref{prop:coiso-sections}(v) and defines, by \eqref{eq:family}, a section with coisotropic image; then apply Proposition \ref{prop:case1-leaves} with $k=-1$.
\end{proof}

\begin{remark}
Equation \eqref{eq:evolHJ} is the Hamilton--Jacobi equation associated with the Herglotz variational principle, studied from the analytical point of view (viscosity solutions, Lax--Oleinik semigroups) in \cite{CCJWY20}. Corollary \ref{cor:evolution} shows that it is the geometric content of the Hamilton--Jacobi problem for $X_H$ with coisotropic sections: $N_\gamma$ is the trace of a Legendrian submanifold moved by the contact flow. On the other hand, if $H\circ\gamma=0$ then $k=0$ and, by Lemma \ref{lem:tangency}, $X_H$ is tangent to every Legendrian leaf: each leaf is a solution of the problem considered in the next subsection.
\end{remark}

\subsection{Case 2: the Hamiltonian vector field and Legendrian sections}

Consider now a section $\gamma\colon Q\to T^*Q\times\R$ of $\pi_1$ and $X_H^\gamma=T\pi_1\circ X_H\circ\gamma$.

\begin{theorem}\label{thm:case2}
Assume that $\gamma(Q)$ is Legendrian, i.e., $\gamma=j^1f$ for a function $f\colon Q\to\R$. Then $X_H$ and $X_H^\gamma$ are $\gamma$-related if and only if
\begin{equation}\label{eq:HJ2}
H\circ j^1f=0,\qquad\text{i.e.,}\qquad H\left(q^i,\pd{f}{q^i},f\right)=0 .
\end{equation}
\end{theorem}

\begin{proof}
It follows from Lemmas \ref{lem:related-tangent} and \ref{lem:tangency}(i). In coordinates, writing $\gamma=(q^i,\gamma_i(q),\gamma_z(q))$, the $\gamma$-relatedness is equivalent to
\[
\pd{H}{q^j}+\gamma_j\pd{H}{z}+\pd{H}{p_i}\pd{\gamma_j}{q^i}=0,\qquad \gamma_i\pd{H}{p_i}-H=\pd{H}{p_i}\pd{\gamma_z}{q^i},
\]
and, when $\gamma_i=\partial\gamma_z/\partial q^i$, the first condition says $d(H\circ\gamma)=0$ and the second one $H\circ\gamma=0$.
\end{proof}

\subsection{Case 3: the evolution vector field and coisotropic sections}

Let $\gamma\colon Q\times\R\to T^*Q\times\R$ be a section of $\pi$ and $\mathcal{E}_H^\gamma=T\pi\circ\mathcal{E}_H\circ\gamma$.

\begin{theorem}\label{thm:case3}
Assume that $N_\gamma=\gamma(Q\times\R)$ is coisotropic. The following assertions are equivalent:
\begin{enumerate}[label=(\roman*)]
\item $\mathcal{E}_H$ and $\mathcal{E}_H^\gamma$ are $\gamma$-related;
\item $A_j^\gamma(H\circ\gamma)=0$ for $j=1,\dots,n$;
\item $d(H\circ\gamma)\wedge\vartheta_\gamma=0$, or equivalently $d(H\circ\gamma)=\pd{(H\circ\gamma)}{z}\,\gamma^*\eta_Q$;
\item $H$ is constant on each Legendrian leaf of $N_\gamma$;
\item locally, $H\big(q,\partial_qS(q,c),S(q,c)\big)=k(c)$ for some function $k$, where $S$ is as in Proposition \ref{prop:coiso-sections}.
\end{enumerate}
\end{theorem}

\begin{proof}
The vector fields $\mathcal{E}_H$ and $X_H$ have the same $\partial_q$ and $\partial_p$ components, and $\mathcal{E}_H^\gamma=\pd{H}{p_i}\partial_{q^i}+\gamma_i\pd{H}{p_i}\partial_z$. As in the proof of Theorem \ref{thm:case1}, (i) is equivalent to
\[
\pd{H}{q^j}+\gamma_j\pd{H}{z}+\pd{H}{p_i}A_i^\gamma(\gamma_j)=0,
\]
which, by \eqref{eq:coiso-gamma}, is $A_j^\gamma(H\circ\gamma)=0$. Since the $A_j^\gamma$ span $\ker\vartheta_\gamma$, (ii) and (iii) are equivalent, and they say that $H\circ\gamma$ is constant on the leaves of $\ker\vartheta_\gamma$, which are mapped by $\gamma$ onto the Legendrian leaves of $N_\gamma$; this is (iv). Finally, (iv) and (v) are equivalent by Proposition \ref{prop:coiso-sections}. Alternatively, (i)$\Leftrightarrow$(iv) follows from Lemmas \ref{lem:related-tangent} and \ref{lem:tangency}(ii), since $\mathcal{E}_H$ takes values in $\Hor$ and therefore it is tangent to $N_\gamma$ if and only if it is tangent to $T N_\gamma\cap\Hor$.
\end{proof}

A section $\gamma$ with coisotropic image satisfying the above conditions is called a \emph{solution of the Hamilton--Jacobi problem for $\mathcal{E}_H$}.

\begin{remark}
Equation (iii) coincides with equation (39) of \cite{dLLMB21}, $d_Q(H\circ\gamma)+\gamma_o\,\gamma^*\theta_Q=0$ with $\gamma_o=\partial(H\circ\gamma)/\partial z$. In \cite[Theorem 4]{dLLMB21} it is assumed that the slices $\gamma(Q\times\{z\})$ are Legendrian.
Taken literally, this condition is very restrictive for a section of $\pi$: the slices lie in the hypersurfaces $z=\text{const}$, on which $\eta_Q$ restricts to $-\rho^*\theta_Q$, so they are Legendrian only if $\bar\gamma=0$, i.e., only if $\gamma$ is the zero section. Theorem \ref{thm:case3} shows that the natural assumption is just the coisotropy of $N_\gamma$, the Legendrian submanifolds being the leaves $j^1S_c(Q)$ rather than the slices.
\end{remark}

\subsection{Case 4: the evolution vector field and Legendrian sections}

\begin{theorem}\label{thm:case4}
Let $\gamma=j^1f\colon Q\to T^*Q\times\R$ be a section of $\pi_1$ with Legendrian image and $\mathcal{E}_H^\gamma=T\pi_1\circ\mathcal{E}_H\circ\gamma$. Then $\mathcal{E}_H$ and $\mathcal{E}_H^\gamma$ are $\gamma$-related if and only if
\begin{equation}\label{eq:HJ4}
d(H\circ j^1f)=0,\qquad\text{i.e.,}\qquad H\left(q^i,\pd{f}{q^i},f\right)=E
\end{equation}
for some constant $E\in\R$.
\end{theorem}

\begin{proof}
It is a direct consequence of Lemmas \ref{lem:related-tangent} and \ref{lem:tangency}(ii). In coordinates, the $\gamma$-relatedness is equivalent to $\pd{H}{q^j}+\gamma_j\pd{H}{z}+\pd{H}{p_i}\pd{\gamma_j}{q^i}=0$ and $\pd{H}{p_i}\big(\pd{\gamma_z}{q^i}-\gamma_i\big)=0$; for $\gamma=j^1f$ the second one is trivially satisfied and the first one is $d(H\circ\gamma)=0$.
\end{proof}

\subsection{Geometric solutions}

The previous results lead us to the following definitions, which do not require the solutions to be horizontal (that is, images of sections), in analogy with \cite[Definition 6]{dLMdDV17}.

\begin{definition}\label{def:geometric-solution}
Let $(M,\eta,H)$ be a contact Hamiltonian system.
\begin{enumerate}[label=(\alph*)]
\item A \emph{Legendrian solution} of the Hamilton--Jacobi problem for $\mathcal{E}_H$ (resp.\ for $X_H$) is a Legendrian submanifold $L\subseteq M$ such that $H|_L$ is constant (resp.\ $H|_L=0$).
\item A \emph{coisotropic solution} of the Hamilton--Jacobi problem for $\mathcal{E}_H$ (resp.\ for $X_H$) is a coisotropic submanifold $N\subseteq M$ of dimension $n+1$ transverse to $\Hor$, such that $\mathcal{E}_H$ (resp.\ $X_H$) is tangent to $N$.
\end{enumerate}
A solution is \emph{horizontal} if, for $M=T^*Q\times\R$, it is the image of a section of $\pi_1$ (case (a)) or of $\pi$ (case (b)).
\end{definition}

By Lemma \ref{lem:tangency}, a Legendrian submanifold $L$ is a Legendrian solution for $\mathcal{E}_H$ (resp.\ $X_H$) if and only if $\mathcal{E}_H$ (resp.\ $X_H$) is tangent to $L$; and by Theorem \ref{thm:case3}, $N$ is a coisotropic solution for $\mathcal{E}_H$ if and only if each Legendrian leaf of $N$ is a Legendrian solution for $\mathcal{E}_H$. Table \ref{tab:cases} summarizes the four cases.

\begin{table}[h]
\centering
\renewcommand{\arraystretch}{1.35}
\small
\begin{tabular}{@{}lll@{}}
\toprule
 & Legendrian sections $j^1f\colon Q\to T^*Q\times\R$ & Coisotropic sections $\gamma\colon Q\times\R\to T^*Q\times\R$\\
\midrule
$X_H$ & Case 2: $H(q,\partial_qf,f)=0$ & Case 1: $H(q,\partial_qS,S)=k(c)\,\partial_cS$\\
 & & \quad (evolutionary: $\partial_tS+H(q,\partial_qS,S)=0$)\\
$\mathcal{E}_H$ & Case 4: $H(q,\partial_qf,f)=E$ & Case 3: $H(q,\partial_qS,S)=k(c)$\\
\bottomrule
\end{tabular}
\medskip
\caption{The four contact Hamilton--Jacobi equations, in terms of the generating functions $f(q)$ of a Legendrian section and $S(q,c)$ of the Legendrian leaves of a coisotropic section.}\label{tab:cases}
\end{table}

\begin{remark}[Complete solutions]
A complete solution for $\mathcal{E}_H$ in the sense of \cite[Definition 6]{dLLMB21} is a local diffeomorphism $\Phi\colon Q\times\R^{n+1}\to T^*Q\times\R$ such that each $\Phi_\lambda=\Phi(\cdot,\lambda)$ is a Legendrian solution. In the present language, a complete solution is a Legendrian foliation of (an open subset of) $T^*Q\times\R$ by solutions, and the functions $f_\alpha=\pi_\alpha\circ\Phi^{-1}$ are first integrals of $\mathcal{E}_H$. Similarly, by Theorem \ref{thm:case3}, any one-parameter subfamily of a complete solution whose union is an $(n+1)$-dimensional submanifold transverse to $\Hor$ is a coisotropic solution for $\mathcal{E}_H$. The reduction results of Section \ref{sec:HJ-reduction} apply to complete solutions which are invariant under the symmetry group, by applying them to each member of the family.
\end{remark}

\section{Symmetries and momentum maps}\label{sec:symmetries}

Many contact Hamiltonian systems arising in mechanics and thermodynamics possess continuous symmetries. As in the symplectic setting, these symmetries give rise to momentum maps, which are the starting point of reduction procedures. In this section we study the properties of these momentum maps that are relevant for the Hamilton--Jacobi theory.

\subsection{Contact actions and momentum maps}

Let $G$ be a connected Lie group with Lie algebra $\g$, acting on the contact manifold $(M,\eta)$ by a left action $\Phi\colon G\times M\to M$ of contactomorphisms, $\Phi_g^*\eta=\eta$. For $\xi\in\g$, let $\xi_M(x)=\frac{d}{dt}\big|_{t=0}\exp(t\xi)\cdot x$ be the infinitesimal generator. We denote by $\Ad^*_g=(\Ad_g)^*\colon\g^*\to\g^*$ the transpose of the adjoint map, so that the coadjoint action is $g\cdot\mu=\Ad^*_{g^{-1}}\mu$, and $\langle\ad^*_\xi\mu,\zeta\rangle=\langle\mu,[\xi,\zeta]\rangle$.

\begin{definition}
The \emph{momentum map} of the action is $J\colon M\to\g^*$ defined by
\[
\langle J(x),\xi\rangle=J_\xi(x)=-\eta_x\big(\xi_M(x)\big),\qquad x\in M,\ \xi\in\g .
\]
\end{definition}

\begin{proposition}\label{prop:momentum}
The momentum map satisfies:
\begin{enumerate}[label=(\roman*)]
\item $dJ_\xi=i_{\xi_M}d\eta$ (the momentum condition);
\item $\xi_M=X_{J_\xi}$, and $\sharpL(dJ_\xi)=\mathcal{E}_{J_\xi}=\xi_M+J_\xi\Reeb$;
\item $[\Reeb,\xi_M]=0$ and $\Reeb(J_\xi)=0$; in particular, the level sets of $J$ are invariant under the flow of $\Reeb$;
\item $J$ is equivariant: $J\circ\Phi_g=\Ad^*_{g^{-1}}\circ J$ for all $g\in G$.
\end{enumerate}
\end{proposition}

\begin{proof}
(i) Since $\Lie_{\xi_M}\eta=0$, $0=d(\eta(\xi_M))+i_{\xi_M}d\eta=-dJ_\xi+i_{\xi_M}d\eta$. (iii) The contactomorphisms $\Phi_g$ preserve $\Reeb$, so $[\xi_M,\Reeb]=0$, and $\Reeb(J_\xi)=-(\Lie_\Reeb\eta)(\xi_M)-\eta([\Reeb,\xi_M])=0$ because $\Lie_\Reeb\eta=0$. (ii) By (i) and (iii), $\flat(\xi_M)=dJ_\xi-J_\xi\eta=dJ_\xi-(\Reeb(J_\xi)+J_\xi)\eta=\flat(X_{J_\xi})$. The second assertion follows from Lemma \ref{lem:sharpL}. (iv) Since $\exp(t\xi)gx=g\exp(t\Ad_{g^{-1}}\xi)x$, we have $\xi_M(gx)=T\Phi_g\big((\Ad_{g^{-1}}\xi)_M(x)\big)$, and then $J_\xi(gx)=-(\Phi_g^*\eta)_x\big((\Ad_{g^{-1}}\xi)_M(x)\big)=\langle J(x),\Ad_{g^{-1}}\xi\rangle$.
\end{proof}

\begin{example}[Cotangent lifted actions]\label{ex:cotangent-lift}
Let $G$ act on $Q$ and consider the cotangent lift $\Phi^{T^*Q}_g(\alpha_q)=\alpha_q\circ T_{gq}\Phi_{g^{-1}}$. The action on $M=T^*Q\times\R$ given by
\[
\Phi_g(\alpha_q,z)=\big(\Phi^{T^*Q}_g(\alpha_q),z\big)
\]
is an action by contactomorphisms, since $(\Phi_g^{T^*Q})^*\theta_Q=\theta_Q$ and therefore $\Phi_g^*\eta_Q=dz-\rho^*\theta_Q=\eta_Q$. As the action is trivial on $z$, $\xi_M=(\xi_{T^*Q},0)$ and
\[
J(\alpha_q,z)=\theta_Q(\xi_{T^*Q})(\alpha_q)=\alpha_q\big(\xi_Q(q)\big)=\langle J^{T^*Q}(\alpha_q),\xi\rangle,
\]
where $J^{T^*Q}\colon T^*Q\to\g^*$ is the standard momentum map of the cotangent lifted action \cite{AM78,MMOPR07}. Hence $J=J^{T^*Q}\circ\rho$ does not depend on $z$ and $J^{-1}(\mu)=(J^{T^*Q})^{-1}(\mu)\times\R$.
\end{example}

\subsection{Regular values}

\begin{proposition}\label{prop:regular}
An element $\mu\in J(M)$ is a regular value of $J$ if and only if, for every $x\in J^{-1}(\mu)$,
\[
\{\xi\in\g\mid \xi_M(x)=-\mu(\xi)\,\Reeb_x\}=\{0\}.
\]
In particular, $0$ is a regular value of $J$ if and only if $\g_x=\{\xi\in\g\mid\xi_M(x)=0\}=\{0\}$ for every $x\in J^{-1}(0)$.
\end{proposition}

\begin{proof}
$T_xJ$ is surjective if and only if the annihilator of its image is zero. Now $\xi$ annihilates $\Im T_xJ$ if and only if $dJ_\xi(x)=i_{\xi_M(x)}d\eta_x=0$, that is, $\xi_M(x)\in\ker d\eta_x=\langle\Reeb_x\rangle$. If $\xi_M(x)=a\Reeb_x$, then $a=\eta_x(\xi_M(x))=-J_\xi(x)=-\mu(\xi)$. Conversely, if $\xi_M(x)=-\mu(\xi)\Reeb_x$, then $dJ_\xi(x)=-\mu(\xi)\,i_{\Reeb}d\eta=0$.
\end{proof}

\begin{remark}
In the symplectic case, $\mu$ is a regular value if and only if the action is locally free on $J^{-1}(\mu)$ \cite[Proposition 1]{dLMdDV17}. In the contact case the infinitesimal generators may be proportional to the Reeb vector field, and this is the reason why the zero value plays a special role. For a cotangent lifted action as in Example \ref{ex:cotangent-lift} which is free on $Q$, every $\mu$ is a regular value: indeed, $\xi_M=(\xi_{T^*Q},0)$ has no $\partial_z$ component, so $\xi_M(x)=-\mu(\xi)\Reeb_x$ forces $\xi_{T^*Q}(\alpha_q)=0$, hence $\xi_Q(q)=0$ and $\xi=0$.
\end{remark}

\subsection{Dissipated momenta}

\begin{proposition}\label{prop:dissipated}
Let $H\in C^\infty(M)$ be $G$-invariant. Then $X_H$ and $\mathcal{E}_H$ are $G$-invariant vector fields and, for every $\xi\in\g$,
\[
X_H(J_\xi)=\mathcal{E}_H(J_\xi)=-\Reeb(H)\,J_\xi .
\]
Consequently, along any integral curve $c(t)$ of $X_H$ or of $\mathcal{E}_H$,
\[
J\big(c(t)\big)=\exp\left(-\int_0^t\Reeb(H)\big(c(s)\big)\,ds\right)J\big(c(0)\big).
\]
In particular, $J^{-1}(0)$ is invariant under the flows of $X_H$ and $\mathcal{E}_H$, while, for a regular value $\mu\neq0$, the level set $J^{-1}(\mu)$ is invariant under $X_H$ (or $\mathcal{E}_H$) if and only if $\Reeb(H)$ vanishes on it.
\end{proposition}

\begin{proof}
The $G$-invariance of $X_H$ and $\mathcal{E}_H$ follows from $\Phi_g^*\eta=\eta$ and $\Phi_g^*H=H$. Using Proposition \ref{prop:momentum}(i) and \eqref{eq:XH},
\[
X_H(J_\xi)=d\eta(\xi_M,X_H)=-\big(dH-\Reeb(H)\eta\big)(\xi_M)=-\xi_M(H)+\Reeb(H)\,\eta(\xi_M)=-\Reeb(H)J_\xi,
\]
since $\xi_M(H)=0$. As $\Reeb(J_\xi)=0$, also $\mathcal{E}_H(J_\xi)=X_H(J_\xi)$. The rest follows by integrating the linear differential equation $\frac{d}{dt}J(c(t))=-\Reeb(H)(c(t))\,J(c(t))$.
\end{proof}

The components of the momentum map are thus \emph{dissipated quantities} in the sense of \cite{dLL20,GGMRR20}. For instance, for a Hamiltonian with linear dissipation $H=H_0+\lambda z$ with $\Reeb(H_0)=0$, the momentum decays as $J(c(t))=e^{-\lambda t}J(c(0))$. Proposition \ref{prop:dissipated} is the first indication that, in the contact setting, the zero level of the momentum map is the natural stage for reduction. It has the following immediate consequence for the Hamilton--Jacobi problem.

\begin{corollary}\label{cor:no-nonzero}
Let $H$ be $G$-invariant and let $L$ be a connected Legendrian solution of the Hamilton--Jacobi problem for $\mathcal{E}_H$ or for $X_H$ (Definition \ref{def:geometric-solution}). If $J|_L=\mu$ is constant and $\mu\neq0$, then $\Reeb(H)|_L=0$. In particular, if $\Reeb(H)$ vanishes nowhere (for instance, for a cotangent lifted action and $H=H_0+\lambda z$ with $\Reeb(H_0)=0$ and $\lambda\neq0$), the momentum map does not take a constant non-zero value on any Legendrian solution.
\end{corollary}

\begin{proof}
By Lemma \ref{lem:tangency}, $\mathcal{E}_H$ (resp.\ $X_H$) is tangent to $L$, and $J_\xi$ is constant on $L$; hence $\mathcal{E}_H(J_\xi)=0$ (resp.\ $X_H(J_\xi)=0$) on $L$. By Proposition \ref{prop:dissipated}, $\Reeb(H)\,\mu(\xi)=0$ on $L$ for every $\xi\in\g$, and it suffices to choose $\xi$ with $\mu(\xi)\neq0$.
\end{proof}

\subsection{Invariant Legendrian submanifolds}

The next result is the contact counterpart of \cite[Lemma 3]{dLMdDV17}, which states that a Lagrangian submanifold of $T^*Q$ is $G$-invariant if and only if the momentum map is constant on it.

\begin{theorem}\label{thm:invariant-legendrian}
Let $L\subseteq M$ be a connected Legendrian submanifold.
\begin{enumerate}[label=(\roman*)]
\item The infinitesimal generators $\xi_M$, $\xi\in\g$, are tangent to $L$ if and only if $J|_L=0$. In particular, if $L$ is a $G$-invariant Legendrian submanifold, then $L\subseteq J^{-1}(0)$; conversely, if $L$ is closed and $J|_L=0$, then $L$ is $G$-invariant.
\item If $J|_L$ is constant, $J|_L=\mu$, then $\mu$ vanishes on $[\g,\g]$ (equivalently, $\Ad^*_g\mu=\mu$ for all $g\in G$) and the \emph{twisted generators} $\xi_M+\mu(\xi)\Reeb$, $\xi\in\g$, are tangent to $L$. Conversely, if the twisted generators associated with some $\mu\in\g^*$ are tangent to $L$, then $J|_L=\mu$.
\end{enumerate}
\end{theorem}

\begin{proof}
By Proposition \ref{prop:momentum}(ii), $\xi_M=X_{J_\xi}$ and $\xi_M+J_\xi\Reeb=\mathcal{E}_{J_\xi}$. (i) follows from Lemma \ref{lem:tangency}(i) applied to $f=J_\xi$ for each $\xi$; for a closed submanifold, tangency of the generators implies invariance under the connected group $G$. (ii) By Lemma \ref{lem:tangency}(ii), $\xi_M+\mu(\xi)\Reeb$ is tangent to $L$ for all $\xi$. Since $\Reeb$ commutes with the generators and $[\xi_M,\zeta_M]=-[\xi,\zeta]_M$ for a left action, the bracket of two twisted generators is $-[\xi,\zeta]_M$, which must be tangent to $L$. But $[\xi,\zeta]_M+\mu([\xi,\zeta])\Reeb$ is also tangent to $L$, hence so is $\mu([\xi,\zeta])\Reeb$. Since $\eta(\Reeb)=1$ and $TL\subseteq\ker\eta$, we conclude that $\mu([\xi,\zeta])=0$. For connected $G$, this means $\ad^*_\xi\mu=0$ for all $\xi$, that is, $\Ad^*_g\mu=\mu$. Conversely, if $\xi_M+\mu(\xi)\Reeb$ is tangent to $L$, then $0=\eta(\xi_M+\mu(\xi)\Reeb)=-J_\xi+\mu(\xi)$ on $L$.
\end{proof}

Theorem \ref{thm:invariant-legendrian} shows a remarkable difference with the symplectic case: a $G$-invariant Legendrian submanifold always lies in the \emph{zero} level of the momentum map, while $G$-invariant Lagrangian submanifolds may lie in any level set $J^{-1}(\mu)$ with $\mu$ a fixed point of the coadjoint action. Legendrian submanifolds on which $J$ takes a constant non-zero value are not $G$-invariant, but invariant under the twisted infinitesimal action $\xi\mapsto\xi_M+\mu(\xi)\Reeb$. The special role of the zero level of the contact momentum map is also well known in contact topology; see, e.g., \cite{Le03}.

\begin{remark}[The Hamilton--Jacobi theory as a particular case]\label{rem:HJ-particular}
As in \cite[Remark 2]{dLMdDV17}, the Hamilton--Jacobi theory for Legendrian sections is a particular case of Theorem \ref{thm:invariant-legendrian}, or rather of Lemma \ref{lem:tangency}, on which it relies. Indeed, assume first that $\Reeb(H)=0$. Then the flow of $X_H$ defines an $\R$-action by contactomorphisms (Proposition \ref{prop:basic-dyn}) whose momentum map is $-\eta(X_H)=H$. By Theorem \ref{thm:invariant-legendrian}(i), a Legendrian submanifold $L$ is invariant under this flow if and only if $H|_L=0$, which is the Hamilton--Jacobi equation \eqref{eq:HJ2}; by Theorem \ref{thm:invariant-legendrian}(ii), $H|_L=E$ if and only if $L$ is invariant under the flow of the twisted generator $X_H+E\Reeb$, which coincides with $\mathcal{E}_H$ on $L$, and this is \eqref{eq:HJ4}. For a general $H$, the flow of $X_H$ consists of conformal contactomorphisms, but the proof of Lemma \ref{lem:tangency} only uses \eqref{eq:XH} restricted to $TL\subseteq\Hor$, and the same conclusions hold. Therefore, invariant solutions of the contact Hamilton--Jacobi problem are invariant Legendrian submanifolds lying in the appropriate level set of the ``momentum map'' $H$.
\end{remark}

\subsection{The geometry of the level sets of the momentum map}

\begin{proposition}\label{prop:levelsets}
Let $\mu$ be a regular value of $J$ and $C=J^{-1}(\mu)$. For every $x\in C$:
\begin{enumerate}[label=(\roman*)]
\item $T_xC=\{v\in T_xM\mid d\eta_x(\xi_M(x),v)=0\ \forall\xi\in\g\}$; in particular, $\Reeb_x\in T_xC$;
\item $T_xC^{\perpL}=\{\xi_M(x)+\mu(\xi)\Reeb_x\mid\xi\in\g\}$;
\item $C$ is coisotropic if and only if $\mu$ vanishes on $[\g,\g]$, i.e., $G_\mu=G$;
\item if $\mu=0$, then $C=J^{-1}(0)$ is a $G$-invariant coisotropic submanifold, $T_xC^{\perpL}=T_x(G\cdot x)$, and $T_x(G\cdot x)\subseteq T_xC\cap\Hor_x$.
\end{enumerate}
\end{proposition}

\begin{proof}
(i) $T_xC=\ker T_xJ=\bigcap_\xi\ker dJ_\xi(x)$, and $dJ_\xi=i_{\xi_M}d\eta$. (ii) Since $\mu$ is regular, $(T_xC)^\circ=\{dJ_\xi(x)\mid\xi\in\g\}$, and $\sharpL(dJ_\xi)=\xi_M+J_\xi\Reeb$ by Proposition \ref{prop:momentum}(ii). (iii) Since $\Reeb_x\in T_xC$, $C$ is coisotropic if and only if $\xi_M(x)\in T_xC$ for all $\xi$. By equivariance, $T_xJ(\xi_M(x))=\frac{d}{dt}\big|_{t=0}\Ad^*_{\exp(-t\xi)}\mu=-\ad^*_\xi\mu$, which vanishes for all $\xi$ if and only if $\mu([\g,\g])=0$. (iv) follows from (ii), (iii) and the equivariance of $J$, since $J_\xi=-\eta(\xi_M)$ vanishes on $C$.
\end{proof}

\begin{remark}
In the symplectic case, $T_x(J^{-1}(\mu))^{\perp}=T_x(G\cdot x)$ for every regular value $\mu$, and $J^{-1}(\mu)$ is coisotropic if $G_\mu=G$ \cite[Lemma 4]{dLMdDV17} (and only if, by the argument in the proof of Proposition \ref{prop:levelsets}(iii)). In the contact case the coisotropy condition is the same, but for $\mu\neq0$ the characteristic directions are not tangent to the orbits: they are generated by the twisted fields $\xi_M+\mu(\xi)\Reeb$, in agreement with Theorem \ref{thm:invariant-legendrian}(ii).
\end{remark}

\begin{example}
For a cotangent lifted action (Example \ref{ex:cotangent-lift}), $C=(J^{T^*Q})^{-1}(\mu)\times\R$. Directly from Lemma \ref{lem:complement}: identifying $\Hor_{(\alpha,z)}$ with $T_\alpha(T^*Q)$ through $T\rho$, $d\eta_Q$ corresponds to $\omega_Q$ and $T C\cap\Hor$ to $T\big((J^{T^*Q})^{-1}(\mu)\big)$. Hence, more generally, $\tilde C\times\R$ is coisotropic in $(T^*Q\times\R,\eta_Q)$ if and only if $\tilde C$ is coisotropic in $(T^*Q,\omega_Q)$.
\end{example}

\section{Coisotropic reduction in contact geometry}\label{sec:reduction}

We now recall the contact analogue of the coisotropic reduction theorem \cite{We77,dLMdDV17}, which was established in \cite{dLL19}. Theorems \ref{thm:coiso-red} and \ref{thm:legendrian-red} below are essentially contained in \cite{dLL19} and \cite[Theorem 7.1 and Proposition 7.7]{dLIL24}, where coisotropic reduction in contact and related geometries is reviewed; we include short proofs, based on Lemma \ref{lem:complement}, since the description of the characteristic distribution as $\ker\vartheta\cap\ker d\vartheta$ and the linear algebra in the proofs will be used in Section \ref{sec:HJ-reduction}.

\begin{theorem}[Contact coisotropic reduction, {\cite{dLL19}, \cite[Theorem 7.1]{dLIL24}}]\label{thm:coiso-red}
Let $(M,\eta)$ be a contact manifold of dimension $2n+1$ and $C\subseteq M$ a coisotropic submanifold of dimension $m$ such that the Reeb vector field $\Reeb$ is tangent to $C$. Let $\vartheta=i_C^*\eta$. Then:
\begin{enumerate}[label=(\roman*)]
\item The characteristic distribution $D=TC^{\perpL}=(TC\cap\Hor)^{\perpd}$ is a regular distribution of rank $2n+1-m$, and
\[
D=\ker\vartheta\cap\ker d\vartheta=\{v\in TC\mid \vartheta(v)=0,\ i_vd\vartheta=0\}.
\]
Consequently, $D$ is involutive and $\Lie_X\vartheta=0$ for every section $X$ of $D$.
\item Assume that the space of leaves $C_{\mathrm{red}}=C/D$ is a manifold such that the projection $\pi_C\colon C\to C_{\mathrm{red}}$ is a submersion. Then $\dim C_{\mathrm{red}}=2m-2n-1$ and there is a unique contact form $\eta_{\mathrm{red}}$ on $C_{\mathrm{red}}$ such that $\pi_C^*\eta_{\mathrm{red}}=\vartheta$. Its Reeb vector field $\Reeb_{\mathrm{red}}$ is $\pi_C$-related to $\Reeb$.
\end{enumerate}
\end{theorem}

\begin{proof}
(i) Since $\Reeb\in TC$, $C$ is transverse to $\Hor$, $T_xC\cap\Hor_x$ has dimension $m-1$, $T_xC=(T_xC\cap\Hor_x)\oplus\langle\Reeb_x\rangle$, and $D_x=(T_xC\cap\Hor_x)^{\perpd}$ (Lemma \ref{lem:complement}) has dimension $2n-m+1$; by coisotropy, $D_x\subseteq T_xC\cap\Hor_x=\ker\vartheta_x$. For $v\in\ker\vartheta_x$, $i_vd\vartheta=0$ if and only if $d\eta(v,w)=0$ for all $w\in T_xC$, and, as $d\eta(v,\Reeb)=0$, if and only if $v\in(T_xC\cap\Hor_x)^{\perpd}$. This proves the expression of $D$. If $X,Y$ are sections of $D$, then $\vartheta([X,Y])=-d\vartheta(X,Y)=0$ and $i_{[X,Y]}d\vartheta=\Lie_X i_Yd\vartheta-i_Y\Lie_Xd\vartheta=-i_Y\,d\,i_Xd\vartheta=0$; moreover $\Lie_X\vartheta=i_Xd\vartheta+d\,i_X\vartheta=0$.

(ii) By (i), $\vartheta$ is basic with respect to the foliation defined by $D$, so there is a unique $1$-form $\eta_{\mathrm{red}}$ with $\pi_C^*\eta_{\mathrm{red}}=\vartheta$. The dimension of $C_{\mathrm{red}}$ is $m-(2n+1-m)$. Let $u=T\pi_C(v)$ with $v\in TC$, and assume that $\eta_{\mathrm{red}}(u)=0$ and $i_ud\eta_{\mathrm{red}}=0$. Then $\vartheta(v)=0$ and $\pi_C^*(i_ud\eta_{\mathrm{red}})=i_vd\vartheta=0$ on $TC$, so $v\in D$ and $u=0$. Hence the radical $\ker d\eta_{\mathrm{red}}$ intersects $\ker\eta_{\mathrm{red}}$ only at $0$. On the other hand, $u_\Reeb=T_x\pi_C(\Reeb_x)$ satisfies $\eta_{\mathrm{red}}(u_\Reeb)=\vartheta(\Reeb)=1$ and $\pi_C^*(i_{u_\Reeb}d\eta_{\mathrm{red}})=i_{\Reeb}d\vartheta=i_C^*(i_\Reeb d\eta)=0$, so that $u_\Reeb$ lies in the radical. Therefore the radical is the line spanned by $u_\Reeb$, which is transverse to $\ker\eta_{\mathrm{red}}$. If $w\in\ker\eta_{\mathrm{red}}$ is $d\eta_{\mathrm{red}}$-orthogonal to $\ker\eta_{\mathrm{red}}$, it is also orthogonal to $u_\Reeb$, hence it lies in the radical and $w=0$. Thus $d\eta_{\mathrm{red}}$ is non-degenerate on $\ker\eta_{\mathrm{red}}$, $\eta_{\mathrm{red}}$ is a contact form, and $u_\Reeb=\Reeb_{\mathrm{red}}(\pi_C(x))$ by the characterization of the Reeb vector field.
\end{proof}

\begin{remark}
If $\Reeb$ is not tangent to $C$, the form $\vartheta$ is not basic in general: for $X\in D$ one only has $\Lie_X\vartheta=i_Xd\vartheta=a\,\vartheta$, so that only the contact \emph{distribution} $\ker\vartheta$ descends to the quotient, and one obtains a reduced contact structure rather than a reduced contact form; compare with \cite{Wi02,dLIL24}. In the applications below $\Reeb$ will always be tangent to $C$ (Proposition \ref{prop:momentum}(iii)).
\end{remark}

\begin{proposition}\label{prop:coiso-dyn}
Under the assumptions of Theorem \ref{thm:coiso-red}(ii), let $H\in C^\infty(M)$ be such that $H|_C$ is constant along the characteristic leaves, i.e., $H|_C=H_{\mathrm{red}}\circ\pi_C$ for some $H_{\mathrm{red}}\in C^\infty(C_{\mathrm{red}})$. Then $X_H$ and $\mathcal{E}_H$ are tangent to $C$, and their restrictions to $C$ are $\pi_C$-related to $X_{H_{\mathrm{red}}}$ and $\mathcal{E}_{H_{\mathrm{red}}}$, respectively.
\end{proposition}

\begin{proof}
Let $x\in C$ and $\alpha\in(T_xC)^\circ$. Then $\sharpL(\alpha)\in D_x\subseteq T_xC\cap\Hor_x$, so $dH(\sharpL(\alpha))=0$, and, by Lemma \ref{lem:sharpL}(i) and \eqref{eq:EH},
\[
\alpha(\mathcal{E}_H)=d\eta(\sharpL(\alpha),\mathcal{E}_H)=-\big(dH-\Reeb(H)\eta\big)(\sharpL(\alpha))=0 .
\]
Hence $\mathcal{E}_H$ is tangent to $C$, and so is $X_H=\mathcal{E}_H-H\Reeb$. Let $u=T_x\pi_C(X_H(x))$ and $v=X_{H_{\mathrm{red}}}(\pi_C(x))$. We have $\eta_{\mathrm{red}}(u)=\vartheta(X_H)=-H(x)=\eta_{\mathrm{red}}(v)$. For $w\in T_xC$, using that $\Reeb$ is tangent to $C$ and $\pi_C$-related to $\Reeb_{\mathrm{red}}$, so that $\Reeb_{\mathrm{red}}(H_{\mathrm{red}})\circ\pi_C=\Reeb(H)|_C$,
\[
d\eta_{\mathrm{red}}(u,T\pi_C w)=d\eta(X_H,w)=dH(w)-\Reeb(H)\eta(w)=d\eta_{\mathrm{red}}(v,T\pi_Cw).
\]
Since $T_x\pi_C$ is onto, $\flat_{\mathrm{red}}(u)=\flat_{\mathrm{red}}(v)$, and $u=v$. Finally, $\mathcal{E}_H=X_H+H\Reeb$ is $\pi_C$-related to $X_{H_{\mathrm{red}}}+H_{\mathrm{red}}\Reeb_{\mathrm{red}}=\mathcal{E}_{H_{\mathrm{red}}}$.
\end{proof}

The reduction of the dynamics in Theorem \ref{thm:zero-reduction} below is a particular case of Proposition \ref{prop:coiso-dyn}.

\begin{remark}
In Proposition \ref{prop:coiso-dyn}, $H$ only needs to be constant on the characteristic leaves of $C$, with no condition outside $C$. This is what allows us to use it in Section \ref{sec:nonzero} for Hamiltonians that are invariant under a twisted action only on a level set of the momentum map.
\end{remark}

\begin{theorem}[Reduction of Legendrian submanifolds, {\cite[Proposition 7.7]{dLIL24}}]\label{thm:legendrian-red}
Under the assumptions of Theorem \ref{thm:coiso-red}(ii), let $L\subseteq M$ be a Legendrian submanifold such that $L$ and $C$ intersect cleanly, and assume that $\pi_C(L\cap C)$ is an (immersed) submanifold of $C_{\mathrm{red}}$ such that $\pi_C\colon L\cap C\to\pi_C(L\cap C)$ is a submersion. Then $L_{\mathrm{red}}=\pi_C(L\cap C)$ is a Legendrian submanifold of $(C_{\mathrm{red}},\eta_{\mathrm{red}})$.

In particular, if $L\subseteq C$, then $D|_L\subseteq TL$, $L$ is a union of characteristic leaves and, whenever $L/D$ is a manifold, $L_{\mathrm{red}}=L/D$ is Legendrian.
\end{theorem}

\begin{proof}
Let $x\in L\cap C$, $V=\Hor_x$, $W=T_xC\cap\Hor_x$, so that $W^{\perpd}=D_x\subseteq W$, and $\Pi=T_xL$, a Lagrangian subspace of $V$. By the clean intersection assumption, $T_x(L\cap C)=T_xL\cap T_xC=\Pi\cap W$. The map $T_x\pi_C$ restricted to $W$ induces an isomorphism $W/W^{\perpd}\cong\ker(\eta_{\mathrm{red}})_{\pi_C(x)}$ which maps the form induced by $d\eta_x$ onto $d\eta_{\mathrm{red}}$, because $\pi_C^*d\eta_{\mathrm{red}}=d\vartheta$. Now $U=\Pi\cap W+W^{\perpd}$ is a Lagrangian subspace of $V$: it is isotropic since $\Pi$ is isotropic and $W^{\perpd}$ is orthogonal to $W$, and, by the modular law (as $W^{\perpd}\subseteq W$),
\[
U^{\perpd}=(\Pi\cap W)^{\perpd}\cap W=(\Pi+W^{\perpd})\cap W=\Pi\cap W+W^{\perpd}=U .
\]
Hence $U/W^{\perpd}$ is Lagrangian in $W/W^{\perpd}$, that is, $T_{\pi_C(x)}L_{\mathrm{red}}=T_x\pi_C(T_x(L\cap C))$ is a Lagrangian subspace of $\ker\eta_{\mathrm{red}}$. Therefore $L_{\mathrm{red}}$ is an integral submanifold of $\ker\eta_{\mathrm{red}}$ of maximal dimension, i.e., a Legendrian submanifold. If $L\subseteq C$, then $T_xL\subseteq W$ and $D_x=W^{\perpd}\subseteq\Pi^{\perpd}=\Pi=T_xL$.
\end{proof}

\subsection{Contact reduction at zero momentum}

\begin{theorem}\label{thm:zero-reduction}
Let $G$ act freely and properly by contactomorphisms on $(M,\eta)$ with momentum map $J$, and assume that $J^{-1}(0)\neq\emptyset$. Then $0$ is a regular value of $J$, $C=J^{-1}(0)$ is a $G$-invariant coisotropic submanifold tangent to $\Reeb$, its characteristic leaves are the $G$-orbits, and $M_0=J^{-1}(0)/G$ is a contact manifold of dimension $2n+1-2\dim G$ with contact form $\eta_0$ characterized by $\pi_0^*\eta_0=i_C^*\eta$, where $\pi_0\colon J^{-1}(0)\to M_0$ is the projection.

If $H\in C^\infty(M)$ is $G$-invariant, let $H_0\in C^\infty(M_0)$ be defined by $H_0\circ\pi_0=H|_{J^{-1}(0)}$. Then $X_H$ and $\mathcal{E}_H$ are tangent to $J^{-1}(0)$ and they are $\pi_0$-related to $X_{H_0}$ and $\mathcal{E}_{H_0}$, respectively.
\end{theorem}

\begin{proof}
The first part follows from Propositions \ref{prop:regular}, \ref{prop:momentum}(iii) and \ref{prop:levelsets}(iv), and from Theorem \ref{thm:coiso-red}: the leaves of $D=T(G\cdot x)$ are the orbits (recall that $G$ is connected), and $J^{-1}(0)/G$ is a manifold because the action is free and proper. The dimension is $2(2n+1-r)-2n-1$ with $r=\dim G$.

Let $H$ be $G$-invariant. By Proposition \ref{prop:dissipated}, $X_H$ is tangent to $C$, and, being $G$-invariant, its restriction to $C$ projects to a vector field $Y=T\pi_0\circ X_H$ on $M_0$. Since $\pi_0^*$ is injective on forms, it suffices to check that the pull-backs of the two sides of \eqref{eq:XH} for $(Y,H_0)$ coincide. We have $\pi_0^*\big(\eta_0(Y)\big)=\eta(X_H)|_C=-H|_C=-\pi_0^*H_0$. Moreover, $\Reeb_0(H_0)\circ\pi_0=\Reeb(H)|_C$ since $\Reeb$ is $\pi_0$-related to $\Reeb_0$, and
\[
\pi_0^*(i_Yd\eta_0)=i_{X_H}\big(i_C^*d\eta\big)=i_C^*\big(dH-\Reeb(H)\eta\big)=\pi_0^*\big(dH_0-\Reeb_0(H_0)\eta_0\big).
\]
Hence $Y=X_{H_0}$. Finally, $\mathcal{E}_H=X_H+H\Reeb$ is $\pi_0$-related to $X_{H_0}+H_0\Reeb_0=\mathcal{E}_{H_0}$.
\end{proof}

This is the contact reduction at zero of Albert, Geiges, Loose and Willett \cite{Al89,Ge97,Lo01,Wi02}, here obtained as a particular case of coisotropic reduction, together with the reduction of the contact dynamics.

\begin{remark}[Reduction at $\mu\neq0$]\label{rem:mu-nonzero}
Let $\mu\neq0$ be a regular value of $J$ with $\mu([\g,\g])=0$. Then $J^{-1}(\mu)$ is coisotropic and tangent to $\Reeb$, and Theorem \ref{thm:coiso-red} applies: its characteristic foliation is generated by the twisted fields $\xi_M+\mu(\xi)\Reeb$. If $\mu$ integrates to a group homomorphism $\chi\colon G\to\R$ (e.g., if $G$ is simply connected), these fields generate the \emph{twisted action} $g\ast x=\mathrm{Fl}^{\Reeb}_{\chi(g)}(\Phi_g(x))$, which is again an action by contactomorphisms, with momentum map $J-\mu$. Thus reduction at $\mu$ for $\Phi$ is reduction at $0$ for the twisted action. However, a $G$-invariant Hamiltonian is invariant under the twisted action only if $\Reeb(H)=0$, and by Proposition \ref{prop:dissipated} $J^{-1}(\mu)$ is not preserved by the dissipative dynamics. This explains why, in the presence of dissipation, reduction at a single value $\mu\neq0$ is not available; instead, one can reduce along the ray $\R^+\mu$. Both procedures are developed for the Hamilton--Jacobi problem in Section \ref{sec:nonzero}. Other approaches to contact reduction at $\mu\neq0$, independent of the contact form, are discussed in \cite{Wi02}.
\end{remark}

\subsection{The cotangent bundle case}

Assume now that $G$ acts freely and properly on $Q$, so that $\pi_G\colon Q\to Q/G$ is a principal $G$-bundle, and consider the lifted action on $M=T^*Q\times\R$ of Example \ref{ex:cotangent-lift}. Recall that $J^{-1}(0)=(J^{T^*Q})^{-1}(0)\times\R$ and that $\alpha_q\in(J^{T^*Q})^{-1}(0)$ if and only if $\alpha_q$ annihilates the vertical space $\ker T_q\pi_G$, that is, if and only if $\alpha_q=(T_q\pi_G)^*\beta$ for a unique $\beta\in T^*_{\pi_G(q)}(Q/G)$.

\begin{proposition}\label{prop:cotangent-red}
The map
\[
\psi_0\colon J^{-1}(0)\to T^*(Q/G)\times\R,\qquad \psi_0\big((T_q\pi_G)^*\beta,z\big)=(\beta,z),
\]
is a surjective submersion whose fibres are the $G$-orbits, and $\psi_0^*\eta_{Q/G}=i^*\eta_Q$, where $i\colon J^{-1}(0)\to T^*Q\times\R$ is the inclusion and $\eta_{Q/G}=dz-\theta_{Q/G}$. Consequently, the reduced contact manifold $(J^{-1}(0)/G,\eta_0)$ is contactomorphic to $(T^*(Q/G)\times\R,\eta_{Q/G})$. If $H$ is a $G$-invariant Hamiltonian on $T^*Q\times\R$, the reduced Hamiltonian is
\[
\widetilde H\colon T^*(Q/G)\times\R\to\R,\qquad \widetilde H(\beta,z)=H\big((T_q\pi_G)^*\beta,z\big),\quad q\in\pi_G^{-1}\big(\pi_{Q/G}(\beta)\big),
\]
and $T\psi_0\circ X_H=X_{\widetilde H}\circ\psi_0$, $T\psi_0\circ\mathcal{E}_H=\mathcal{E}_{\widetilde H}\circ\psi_0$ on $J^{-1}(0)$.
\end{proposition}

\begin{proof}
If $\alpha_q=\beta\circ T_q\pi_G$, then $\Phi_g^{T^*Q}(\alpha_q)=\beta\circ T_q\pi_G\circ T_{gq}\Phi_{g^{-1}}=\beta\circ T_{gq}\pi_G$, so $\psi_0$ is constant on orbits; conversely, if $\psi_0(\alpha_q,z)=\psi_0(\alpha'_{q'},z')$ then $z=z'$ and $q'=gq$ for some $g$, and both $\Phi^{T^*Q}_g(\alpha_q)$ and $\alpha'_{q'}$ are equal to $\beta\circ T_{q'}\pi_G$. That $\psi_0$ is a surjective submersion is the well-known cotangent bundle reduction at zero \cite{AM78,MMOPR07}, see also \cite[Appendix B]{dLMdDV17}. For $v\in T_{(\alpha_q,z)}J^{-1}(0)$, denoting by $\pi_1$ and $\tilde\pi_1$ the projections onto $Q$ and $Q/G$, and using $\tilde\pi_1\circ\psi_0=\pi_G\circ\pi_1$,
\[
(\psi_0^*\theta_{Q/G})(v)=\beta\big(T\tilde\pi_1(T\psi_0 v)\big)=\beta\big(T\pi_G(T\pi_1 v)\big)=\alpha_q(T\pi_1v)=(\rho^*\theta_Q)(v),
\]
and $\psi_0^*dz=dz$. Hence $\psi_0^*\eta_{Q/G}=i^*\eta_Q$, and $\psi_0$ induces a diffeomorphism $J^{-1}(0)/G\to T^*(Q/G)\times\R$ which is a contactomorphism by uniqueness of $\eta_0$. The statements about the Hamiltonian and the dynamics follow from Theorem \ref{thm:zero-reduction}. (Compare with the cosphere bundle reduction of \cite{DOR03}.)
\end{proof}

\begin{remark}
In \cite{dLMdDV17}, the reduced spaces $(J^{T^*Q})^{-1}(\mu)/G$ are symplectomorphic to $T^*(Q/G)$ endowed with the canonical symplectic form plus a magnetic term $B_\mu$, which vanishes for $\mu=0$. In the contact setting, $G$-invariant solutions only occur at $\mu=0$, where the reduced space is the canonical contact manifold $T^*(Q/G)\times\R$; non-zero momenta are treated in Section \ref{sec:nonzero}, where again no magnetic terms appear.
\end{remark}

\section{Reduction and reconstruction of the contact Hamilton--Jacobi equation}\label{sec:HJ-reduction}

Throughout this section $G$ acts freely and properly on $Q$, the action is lifted to $T^*Q\times\R$ as in Example \ref{ex:cotangent-lift}, $H\colon T^*Q\times\R\to\R$ is a $G$-invariant Hamiltonian and $\widetilde H\colon T^*(Q/G)\times\R\to\R$ is the reduced Hamiltonian of Proposition \ref{prop:cotangent-red}. We denote $r=\dim G$. The situation is summarized in the diagram
\[
\begin{tikzcd}
T^*Q\times\R \arrow[d,"\pi_1"'] & J^{-1}(0) \arrow[l,hook'] \arrow[r,"\psi_0"] & T^*(Q/G)\times\R \arrow[d,"\tilde\pi_1"]\\
Q \arrow[rr,"\pi_G"] & & Q/G
\end{tikzcd}
\]

\subsection{Reduction of Legendrian solutions (Cases 2 and 4)}

\begin{theorem}[Reduction]\label{thm:HJ-red}
Let $L\subseteq T^*Q\times\R$ be a closed $G$-invariant Legendrian solution of the Hamilton--Jacobi problem for $\mathcal{E}_H$, with $H|_L=E$ (resp.\ a solution for $X_H$, i.e., $E=0$). Then:
\begin{enumerate}[label=(\roman*)]
\item $L\subseteq J^{-1}(0)$;
\item $\widetilde L=\psi_0(L)$ is a Legendrian submanifold of $(T^*(Q/G)\times\R,\eta_{Q/G})$ and a Legendrian solution of the Hamilton--Jacobi problem for $\mathcal{E}_{\widetilde H}$ (resp.\ $X_{\widetilde H}$), with $\widetilde H|_{\widetilde L}=E$;
\item if $L=j^1S(Q)$ is horizontal, then $S$ is $G$-invariant, $S=\widetilde S\circ\pi_G$ for a function $\widetilde S\colon Q/G\to\R$, and $\widetilde L=j^1\widetilde S(Q/G)$ is horizontal. The reduced Hamilton--Jacobi equation is
\[
\widetilde H\big(j^1\widetilde S\big)=E\quad(\text{resp. } =0)\qquad\text{on } Q/G .
\]
\end{enumerate}
\end{theorem}

\begin{proof}
(i) is Theorem \ref{thm:invariant-legendrian}(i). (ii) Since $L\subseteq J^{-1}(0)$ is $G$-invariant and $G$ acts freely and properly, $\widetilde L=\psi_0(L)\cong L/G$ is a submanifold of dimension $n-r=\dim(Q/G)$, and $\psi_0|_L\colon L\to\widetilde L$ is a submersion. As $\psi_0^*\eta_{Q/G}=i^*\eta_Q$ vanishes on $TL$, $\eta_{Q/G}$ vanishes on $T\widetilde L=T\psi_0(TL)$. Thus $\widetilde L$ is Legendrian (alternatively, apply Theorem \ref{thm:legendrian-red} with $L\subseteq C$). Moreover $\widetilde H\circ\psi_0=H$ on $J^{-1}(0)$, so $\widetilde H|_{\widetilde L}=E$. (iii) If $L=j^1S(Q)$ is $G$-invariant, then $\Phi_g(j^1S(q))=\big(\Phi^{T^*Q}_g(dS(q)),S(q)\big)\in L$ must be equal to $j^1S(gq)$, so $S(gq)=S(q)$. Since $\pi_G$ is a surjective submersion with connected fibres, $S=\widetilde S\circ\pi_G$ for a smooth $\widetilde S$, and $dS(q)=(T_q\pi_G)^*d\widetilde S(\pi_G(q))$. Hence $\psi_0(j^1S(q))=j^1\widetilde S(\pi_G(q))$.
\end{proof}

\begin{theorem}[Reconstruction]\label{thm:HJ-rec}
Let $\widetilde L\subseteq T^*(Q/G)\times\R$ be a Legendrian solution of the Hamilton--Jacobi problem for $\mathcal{E}_{\widetilde H}$ with $\widetilde H|_{\widetilde L}=E$ (resp.\ for $X_{\widetilde H}$). Then $L=\psi_0^{-1}(\widetilde L)\subseteq J^{-1}(0)$ is a $G$-invariant Legendrian solution of the Hamilton--Jacobi problem for $\mathcal{E}_H$ with $H|_L=E$ (resp.\ for $X_H$). If $\widetilde L=j^1\widetilde S(Q/G)$, then $L=j^1(\widetilde S\circ\pi_G)(Q)$.
\end{theorem}

\begin{proof}
Since $\psi_0$ is a submersion, $L$ is a submanifold of dimension $(n-r)+r=n$, and it is $G$-invariant because the fibres of $\psi_0$ are the $G$-orbits. As $i^*\eta_Q=\psi_0^*\eta_{Q/G}$ and $\eta_{Q/G}$ vanishes on $T\widetilde L$, $\eta_Q$ vanishes on $TL$; hence $L$ is Legendrian. Finally $H|_L=\widetilde H\circ\psi_0|_L=E$. The last assertion follows as in the proof of Theorem \ref{thm:HJ-red}(iii), since $j^1(\widetilde S\circ\pi_G)$ takes values in $J^{-1}(0)$.
\end{proof}

\begin{corollary}\label{cor:bijection}
There is a one-to-one correspondence
\[
\left\{\begin{array}{c}\text{closed $G$-invariant}\\ \text{Legendrian solutions of the}\\ \text{HJ problem for } \mathcal{E}_H \ (\text{resp. } X_H)\end{array}\right\}
\ \longleftrightarrow\
\left\{\begin{array}{c}\text{closed Legendrian}\\ \text{solutions of the}\\ \text{HJ problem for } \mathcal{E}_{\widetilde H}\ (\text{resp. } X_{\widetilde H})\end{array}\right\},
\]
$L\mapsto\psi_0(L)$, which preserves the value $E$ of the Hamiltonian and horizontality. For horizontal solutions it reads $S=\widetilde S\circ\pi_G$.
\end{corollary}

\begin{remark}
In the symplectic case \cite[Remark 6]{dLMdDV17}, a $G$-invariant horizontal Lagrangian submanifold $\Im dS$ has a $G$-invariant generating function only when $J(\Im dS)=0$; for $\mu\neq0$ one obtains invariant Lagrangian submanifolds with non-invariant generating functions, and the reduced equation involves a magnetic term, becoming an algebraic-PDE system. In the contact case, by Theorem \ref{thm:invariant-legendrian}, $G$-invariant Legendrian submanifolds automatically lie in $J^{-1}(0)$ and the generating function of a horizontal one is always $G$-invariant: the reduced Hamilton--Jacobi equation is again a single first order PDE on $Q/G$, with $\dim G$ fewer independent variables.
\end{remark}

\subsection{Reduction of coisotropic solutions (Cases 1 and 3)}

Let $\gamma\colon Q\times\R\to T^*Q\times\R$ be a section of $\pi$. We say that $\gamma$ is \emph{$G$-equivariant} if $\gamma(gq,z)=\Phi_g(\gamma(q,z))$, where $G$ acts on $Q\times\R$ by $g(q,z)=(gq,z)$.

\begin{theorem}\label{thm:HJ-red-coiso}
Let $\gamma$ be a $G$-equivariant section of $\pi$ with coisotropic image taking values in $J^{-1}(0)$. Then:
\begin{enumerate}[label=(\roman*)]
\item each Legendrian leaf of $N_\gamma=\gamma(Q\times\R)$ is (locally) $G$-invariant; equivalently, the functions $S(q,c)$ of Proposition \ref{prop:coiso-sections} can be chosen $G$-invariant in $q$;
\item there is a unique section $\tilde\gamma\colon(Q/G)\times\R\to T^*(Q/G)\times\R$ such that $\tilde\gamma\circ(\pi_G\times\id_\R)=\psi_0\circ\gamma$, and $\tilde\gamma$ has coisotropic image;
\item $\gamma$ is a solution of the Hamilton--Jacobi problem for $X_H$ (resp.\ $\mathcal{E}_H$) if and only if $\tilde\gamma$ is a solution of the Hamilton--Jacobi problem for $X_{\widetilde H}$ (resp.\ $\mathcal{E}_{\widetilde H}$).
\end{enumerate}
Conversely, every section $\tilde\gamma$ with coisotropic image solving the reduced problem is obtained in this way from a unique $G$-equivariant section $\gamma$ with values in $J^{-1}(0)$, namely
\[
\gamma(q,z)=\big((T_q\pi_G)^*\tilde\gamma_{Q/G}(\pi_G(q),z),\,z\big),\qquad\text{where } \tilde\gamma=(\tilde\gamma_{Q/G},z).
\]
\end{theorem}

\begin{proof}
(i) For $\xi\in\g$, $\vartheta_\gamma(\xi_Q,0)=\eta_Q\big(T\gamma(\xi_Q,0)\big)=\eta_Q(\xi_M\circ\gamma)=-J_\xi\circ\gamma=0$, using the equivariance of $\gamma$. Hence the generators of the action on $Q\times\R$ are tangent to the leaves of $\ker\vartheta_\gamma$, and $\gamma$ maps them onto the Legendrian leaves. (ii) The map $\psi_0\circ\gamma$ is $G$-invariant, so it factors through the submersion $\pi_G\times\id_\R$; $\tilde\gamma$ is a section since $\psi_0$ preserves $z$ and $\tilde\pi_1\circ\psi_0=\pi_G\circ\pi_1$. Moreover
\[
(\pi_G\times\id_\R)^*\vartheta_{\tilde\gamma}=\gamma^*\psi_0^*\eta_{Q/G}=\gamma^*\eta_Q=\vartheta_\gamma,
\]
so $(\pi_G\times\id_\R)^*(\vartheta_{\tilde\gamma}\wedge d\vartheta_{\tilde\gamma})=\vartheta_\gamma\wedge d\vartheta_\gamma=0$, and, $\pi_G\times\id_\R$ being a surjective submersion, $\vartheta_{\tilde\gamma}\wedge d\vartheta_{\tilde\gamma}=0$; by Proposition \ref{prop:coiso-sections}, $\tilde\gamma$ has coisotropic image. (iii) We have $H\circ\gamma=\widetilde H\circ\tilde\gamma\circ(\pi_G\times\id_\R)$, and the vector field $\partial_z$ on $Q\times\R$ is $(\pi_G\times\id_\R)$-related to $\partial_z$ on $(Q/G)\times\R$. Hence the left-hand sides of \eqref{eq:HJ1-intr} and of Theorem \ref{thm:case3}(iii) for $\gamma$ are the pull-backs by $\pi_G\times\id_\R$ of the corresponding forms for $\tilde\gamma$, and the equivalence follows from the injectivity of the pull-back by a surjective submersion. The converse is proved as Theorem \ref{thm:HJ-rec}.
\end{proof}

\begin{remark}
In terms of generating functions, Theorem \ref{thm:HJ-red-coiso} says that $G$-invariant families $S(q,c)=\widetilde S(\pi_G(q),c)$ solve $H(q,\partial_qS,S)=k(c)\partial_cS$ (resp.\ $=k(c)$) if and only if $\widetilde S$ solves the same equation for $\widetilde H$. For Case 1 this includes the reduction of the evolutionary equation \eqref{eq:evolHJ}: $G$-invariant solutions $S(q,t)$ of $\partial_tS+H(q,\partial_qS,S)=0$ correspond to solutions of $\partial_t\widetilde S+\widetilde H(x,\partial_x\widetilde S,\widetilde S)=0$ on $Q/G$.
\end{remark}

\subsection{Reduction of the projected dynamics and reconstruction}

As in \cite[Section 4]{dLMdDV17}, the reduction of the Hamilton--Jacobi equation is compatible with the reduction of the projected dynamics. We state it for Case 4; the other cases are analogous.

\begin{proposition}\label{prop:reconstruction}
Let $\gamma=j^1S$ be a $G$-invariant horizontal Legendrian solution for $\mathcal{E}_H$, and $\tilde\gamma=j^1\widetilde S$ its reduction. Then the projected vector field $Y=\mathcal{E}_H^\gamma=T\pi_1\circ\mathcal{E}_H\circ\gamma\in\X(Q)$ is $G$-invariant and $\pi_G$-related to $\widetilde Y=\mathcal{E}_{\widetilde H}^{\tilde\gamma}\in\X(Q/G)$. Moreover, let $\mathcal{A}\colon TQ\to\g$ be a principal connection on $\pi_G\colon Q\to Q/G$, let $\tilde q(t)$ be the integral curve of $\widetilde Y$ through $\pi_G(q_0)$ and $q_h(t)$ its horizontal lift with $q_h(0)=q_0$. If $g(t)$ is the solution of
\[
\dot g(t)=T_eL_{g(t)}\Big(\mathcal{A}\big(Y(q_h(t))\big)\Big),\qquad g(0)=e,
\]
then $q(t)=g(t)\cdot q_h(t)$ is the integral curve of $Y$ through $q_0$, and $\gamma(q(t))$ is the integral curve of $\mathcal{E}_H$ through $j^1S(q_0)$.
\end{proposition}

\begin{proof}
The $G$-invariance of $Y$ follows from the invariance of $\mathcal{E}_H$ and of $\gamma$. Using $\pi_G\circ\pi_1=\tilde\pi_1\circ\psi_0$ on $J^{-1}(0)$, Proposition \ref{prop:cotangent-red} and Theorem \ref{thm:HJ-red}(iii),
\[
T\pi_G(Y(q))=T\tilde\pi_1\,T\psi_0\big(\mathcal{E}_H(\gamma(q))\big)=T\tilde\pi_1\,\mathcal{E}_{\widetilde H}\big(\tilde\gamma(\pi_G(q))\big)=\widetilde Y(\pi_G(q)).
\]
For the reconstruction, write $q(t)=g(t)q_h(t)$. Then
\[
T\Phi_{g^{-1}}\dot q=\dot q_h+\big(g^{-1}\dot g\big)_Q(q_h),
\]
while $T\Phi_{g^{-1}}Y(q)=Y(q_h)$ by invariance. Since both $Y(q_h(t))$ and $\dot q_h(t)$ project onto $\dot{\tilde q}(t)$, their difference is vertical, and as $\mathcal{A}(\dot q_h)=0$ and $\mathcal{A}(\xi_Q)=\xi$, the equation $\dot q=Y(q)$ is equivalent to $g^{-1}\dot g=\mathcal{A}(Y(q_h))$. The last assertion is the Hamilton--Jacobi property of $\gamma$ (Theorem \ref{thm:case4}).
\end{proof}

The situation is summarized in the following diagram, analogous to \cite[Fig.~1]{dLMdDV17}:
\[
\begin{tikzcd}[column sep=8em,row sep=3.5em]
(T^*Q\times\R,\ \mathcal{E}_H) \arrow[r,"\text{reduction}"] \arrow[d,"\text{Hamilton--Jacobi}"'] & (T^*(Q/G)\times\R,\ \mathcal{E}_{\widetilde H}) \arrow[d,"\text{Hamilton--Jacobi}"] \arrow[l,bend right=15,dashed,"\text{reconstruction}"']\\
(Q,\ \mathcal{E}^{\gamma}_H) \arrow[r,"\pi_G"] & (Q/G,\ \mathcal{E}^{\tilde\gamma}_{\widetilde H}) \arrow[l,bend left=15,dashed,"\text{reconstruction}"]
\end{tikzcd}
\]

\begin{remark}\label{rem:mech-conn}
The same reconstruction scheme, based on a principal connection, is used for contact Lagrangian systems in \cite{ACdLSS24}. For natural Hamiltonians the reconstruction at zero momentum is trivial. Indeed, assume that $H=\frac12\|p\|^2+U(q,z)$, where $\|\cdot\|$ is the norm defined by a $G$-invariant Riemannian metric $\langle\cdot,\cdot\rangle$ on $Q$, and let $\mathcal{A}$ be the mechanical connection, whose horizontal distribution is the orthogonal complement of the orbits. Then $Y=\operatorname{grad}S$, and $\langle Y,\xi_Q\rangle=dS(\xi_Q)=0$ because $S$ is $G$-invariant. Hence $\mathcal{A}(Y)=0$ and $g(t)=e$, as in the conservative case at zero momentum. Non-trivial reconstruction phases require non-natural Hamiltonians, for instance with gyroscopic terms linear in the momenta.
\end{remark}

\subsection{Lagrangian slices and symplectic coisotropic reduction}

In the formulation of \cite{dLLMB21} for the Hamiltonian vector field, each slice $\gamma^z=\rho\circ\gamma(\cdot,z)$ is a closed $1$-form on $Q$. One can then reduce the slices by means of the symplectic coisotropic reduction theorem \cite[Theorem 5]{dLMdDV17}. We include this approach for completeness and to compare it with the contact one.

\begin{proposition}\label{prop:slices}
Let $\gamma$ be a solution of the Hamilton--Jacobi problem for $X_H$ (or $\mathcal{E}_H$) such that the slices $\gamma^z(Q)$ are Lagrangian submanifolds of $(T^*Q,\omega_Q)$.
\begin{enumerate}[label=(\roman*)]
\item Let $\tilde C\subseteq T^*Q$ be a coisotropic submanifold such that $\tilde C/\!\sim$ is a symplectic manifold $(\tilde C/\!\sim,\omega_{\tilde C})$ with projection $\tilde\pi$. If, for $z\in\R$, $\gamma^z(Q)\cap\tilde C$ is a clean intersection, then $\tilde\pi(\gamma^z(Q)\cap\tilde C)$ is a Lagrangian submanifold of $\tilde C/\!\sim$ (whenever it is a submanifold).
\item If $\gamma$ is $G$-equivariant, then $J^{T^*Q}$ is constant on each slice, $J^{T^*Q}(\gamma^z(Q))=\mu(z)$, and $\mu(z)$ vanishes on $[\g,\g]$; hence the slice projects onto a Lagrangian submanifold of $(J^{T^*Q})^{-1}(\mu(z))/G\cong\big(T^*(Q/G),\omega_{Q/G}+B_{\mu(z)}\big)$. Moreover, on any interval $I\ni z_0$ on which $\mu(z)\neq0$, $\bar\gamma$ does not vanish on $Q\times I$, $\partial_z\bar\gamma=\sigma\bar\gamma$ with $\sigma=\sigma(z)$ independent of $q$, and $\mu(z)=\exp\big(\int_{z_0}^z\sigma\big)\mu(z_0)$ lies on the ray $\R^+\mu(z_0)$.
\item If $\gamma$ takes values in $J^{-1}(0)$, i.e., $\mu\equiv0$, these slice-wise reductions are the slices of the section $\tilde\gamma$ of Theorem \ref{thm:HJ-red-coiso}.
\end{enumerate}
\end{proposition}

\begin{proof}
(i) is the symplectic coisotropic reduction theorem \cite{We77}, \cite[Theorem 5]{dLMdDV17}, applied to the Lagrangian submanifold $\gamma^z(Q)$. (ii) Each slice is a $G$-invariant Lagrangian submanifold, so $J^{T^*Q}$ is constant on it by \cite[Lemma 3]{dLMdDV17}; by equivariance, $\mu(z)=\Ad^*_{g^{-1}}\mu(z)$ for all $g\in G$, which for connected $G$ means that $\mu(z)$ vanishes on $[\g,\g]$. The reduction follows from \cite[Section 3.3]{dLMdDV17}. Since $\langle\mu(z),\xi\rangle=\bar\gamma(\xi_Q)(q,z)$ for every $q$, $\bar\gamma$ does not vanish on $Q\times I$. By Remark \ref{rem:lagrangian-slices}, $\partial_z\bar\gamma=\sigma\bar\gamma$ there; hence $\partial_z\langle\mu(z),\xi\rangle=\partial_z\big(\bar\gamma(\xi_Q)\big)=\sigma\langle\mu(z),\xi\rangle$. Since the left-hand side does not depend on $q$, choosing $\xi$ with $\langle\mu(z),\xi\rangle\neq0$ we see that $\sigma$ depends only on $z$, and the claim follows by integration. (iii) is a direct consequence of the definitions.
\end{proof}

\begin{remark}
Proposition \ref{prop:slices} reduces the symplectic slices of a solution but it does not reduce the contact dynamics, which mixes the slices. The values $\mu(z)$ of the momentum map on the slices move along a ray of $\g^*$, which agrees with the exponential decay of the momentum map in Proposition \ref{prop:dissipated}, and is reminiscent of the use of rays $\R^+\mu$ in the contact reduction of Willett \cite{Wi02}. The contact reduction of Theorem \ref{thm:HJ-red-coiso}, which requires $\mu\equiv0$, is the one compatible with the dissipative dynamics at zero momentum; when the momenta move along a ray, the appropriate tool is the reduction along rays of Theorem \ref{thm:HJ-ray}.
\end{remark}

\section{Reduction at non-zero momentum}\label{sec:nonzero}

By Theorem \ref{thm:invariant-legendrian}, Legendrian solutions on which the momentum map takes a constant value $\mu\neq0$ are never $G$-invariant, and by Corollary \ref{cor:no-nonzero} they can only exist where $\Reeb(H)=0$. On the other hand, by Proposition \ref{prop:dissipated}, along the dissipative dynamics the momentum map moves along a ray $\R^+\mu$. In this section we develop the two corresponding reduction procedures for the Hamilton--Jacobi problem: a \emph{twisted reduction} at $\mu$, for Hamiltonians with $\Reeb(H)=0$ on $J^{-1}(\mu)$, and a \emph{reduction along the ray} $\R^+\mu$, valid for every $G$-invariant Hamiltonian.

As in Section \ref{sec:HJ-reduction}, the connected Lie group $G$, of dimension $r$, acts freely and properly on $Q$, and its action is lifted to $M=T^*Q\times\R$ as in Example \ref{ex:cotangent-lift}. We fix $\mu\in\g^*$, $\mu\neq0$, vanishing on $[\g,\g]$, so that $\mu\colon\g\to\R$ is a Lie algebra homomorphism, and we assume:
\begin{itemize}
\item[(A)] there is a Lie group homomorphism $\chi\colon G\to\R$ with $T_e\chi=\mu$.
\end{itemize}
Assumption (A) is automatically satisfied if $G$ is simply connected. It fails for compact $G$, which admits no non-trivial homomorphism into $\R$; for instance, it fails for the action of $SO(2)$ in the central force problem of Section \ref{sec:examples}. There, a Legendrian submanifold with $J=\mu\neq0$ would have, locally, a generating function $S=\mu\,\phi+f(r)$, with $\phi$ the polar angle, which is multivalued, so it cannot be horizontal. Assumption (A) is only used in Section \ref{subsec:twisted} and to guarantee that $K$ is closed; see Remark \ref{rem:ray}(e). We put $\mathfrak k=\ker\mu$, $K=(\ker\chi)_0$ (the identity component of $\ker\chi$), and we fix $\hat\zeta\in\g$ with $\mu(\hat\zeta)=1$.

\begin{lemma}\label{lem:character}
\begin{enumerate}[label=(\roman*)]
\item $\chi$ is surjective, $K$ is a closed connected normal subgroup of $G$ with Lie algebra $\mathfrak k$, and $\chi$ induces an isomorphism of Lie groups $\bar\chi\colon G/K\to\R$.
\item $K$ acts freely and properly on $Q$; denoting by $\pi_K\colon Q\to Q/K$ the projection, $G/K\cong\R$ acts freely and properly on $Q/K$, and $Q/K\to Q/G$ is a principal $\R$-bundle.
\item There is a function $\varphi\in C^\infty(Q)$ such that $\varphi(gq)=\varphi(q)+\chi(g)$ for all $g\in G$, $q\in Q$. For such a function, the $1$-form $d\varphi$ is $G$-invariant and $d\varphi(\xi_Q)=\mu(\xi)$ for all $\xi\in\g$.
\end{enumerate}
\end{lemma}

\begin{proof}
(i) $\chi(G)$ is a connected non-trivial subgroup of $\R$, hence $\chi(G)=\R$. The kernel $\ker\chi$ is a closed normal subgroup with Lie algebra $\mathfrak k$, so $K$ is closed, connected and normal, and $(\ker\chi)/K$ is discrete. Then $G/K$ is a connected one-dimensional Lie group and $\bar\chi\colon G/K\to\R$ is a surjective homomorphism with discrete kernel. Since a compact group cannot be mapped onto $\R$, $G/K\cong\R$, and a surjective homomorphism $\R\to\R$ is a linear isomorphism. (ii) $K$ is a closed subgroup, so it acts freely and properly. If $gK\cdot\pi_K(q)=\pi_K(q)$, then $gq=kq$ for some $k\in K$, hence $g=k$; thus $G/K$ acts freely, and properly because $G$ acts properly. The orbit space of this action is $Q/G$. (iii) A principal bundle with structure group $\R$ is trivial, so there is a $\R$-equivariant function $\bar\varphi\colon Q/K\to\R$, i.e., $\bar\varphi(gK\cdot y)=\bar\varphi(y)+\bar\chi(gK)$; we take $\varphi=\bar\varphi\circ\pi_K$. Differentiating $\varphi\circ\Phi_g=\varphi+\chi(g)$ gives $\Phi_g^*d\varphi=d\varphi$, and differentiating $\varphi(\exp(t\xi)q)=\varphi(q)+t\mu(\xi)$ at $t=0$ gives $d\varphi(\xi_Q)=\mu(\xi)$.
\end{proof}

\subsection{Twisted reduction: the conservative case}\label{subsec:twisted}

Let $G$ act on $M$ by the \emph{twisted action} of Remark \ref{rem:mu-nonzero},
\[
g\star(\alpha_q,z)=\big(\Phi^{T^*Q}_g(\alpha_q),\,z+\chi(g)\big),
\]
which is an action by contactomorphisms with infinitesimal generators $\xi_M+\mu(\xi)\Reeb$ and momentum map $J-\mu$. Consider the map
\[
T_\varphi\colon M\to M,\qquad T_\varphi(\alpha_q,z)=\big(\alpha_q-d\varphi(q),\,z-\varphi(q)\big).
\]

\begin{proposition}\label{prop:shift}
\begin{enumerate}[label=(\roman*)]
\item $T_\varphi$ is a strict contactomorphism, $T_\varphi^*\eta_Q=\eta_Q$, it satisfies $J\circ T_\varphi=J-\mu$ and $T_\varphi(g\star x)=\Phi_g(T_\varphi(x))$, and it maps $J^{-1}(\mu)$ onto $J^{-1}(0)$.
\item The map $\psi_\mu=\psi_0\circ T_\varphi\colon J^{-1}(\mu)\to T^*(Q/G)\times\R$, that is,
\[
\psi_\mu(\alpha_q,z)=\big(\beta,\,z-\varphi(q)\big),\qquad (T_q\pi_G)^*\beta=\alpha_q-d\varphi(q),
\]
is a surjective submersion whose fibres are the orbits of the twisted action, and $\psi_\mu^*\eta_{Q/G}=i_\mu^*\eta_Q$, where $i_\mu\colon J^{-1}(\mu)\to M$ is the inclusion. Consequently, the reduced space of the coisotropic submanifold $J^{-1}(\mu)$ is contactomorphic to $(T^*(Q/G)\times\R,\eta_{Q/G})$.
\item Let $H$ be $G$-invariant. Then $H|_{J^{-1}(\mu)}$ is invariant under the twisted action if and only if $\Reeb(H)=0$ on $J^{-1}(\mu)$. In this case, $H=\widetilde H_\mu\circ\psi_\mu$ on $J^{-1}(\mu)$, where
\[
\widetilde H_\mu(\beta,z)=H\big((T_q\pi_G)^*\beta+d\varphi(q),\,z+\varphi(q)\big),\qquad q\in\pi_G^{-1}\big(\pi_{Q/G}(\beta)\big),
\]
$X_H$ and $\mathcal{E}_H$ are tangent to $J^{-1}(\mu)$ and $\psi_\mu$-related to $X_{\widetilde H_\mu}$ and $\mathcal{E}_{\widetilde H_\mu}$, respectively, and $\Reeb(\widetilde H_\mu)=0$.
\end{enumerate}
\end{proposition}

\begin{proof}
(i) $T_\varphi$ is a diffeomorphism, with inverse $(\alpha_q,z)\mapsto(\alpha_q+d\varphi(q),z+\varphi(q))$. The fibre translation $\tau(\alpha_q)=\alpha_q-d\varphi(q)$ of $T^*Q$ satisfies $\tau^*\theta_Q=\theta_Q-\pi_Q^*d\varphi$, hence
\[
T_\varphi^*\eta_Q=dz-\pi_1^*d\varphi-\rho^*\theta_Q+\pi_1^*d\varphi=\eta_Q .
\]
Moreover $J_\xi(T_\varphi(\alpha_q,z))=\alpha_q(\xi_Q)-d\varphi(\xi_Q)=J_\xi(\alpha_q,z)-\mu(\xi)$, and, since $d\varphi$ is $G$-invariant (so that $\Phi^{T^*Q}_g(d\varphi(q))=d\varphi(gq)$) and $\varphi(gq)=\varphi(q)+\chi(g)$,
\[
T_\varphi(g\star(\alpha_q,z))=\big(\Phi^{T^*Q}_g(\alpha_q)-d\varphi(gq),\,z+\chi(g)-\varphi(gq)\big)=\Phi_g\big(T_\varphi(\alpha_q,z)\big).
\]
(ii) follows from (i) and Proposition \ref{prop:cotangent-red}. The last assertion follows from the uniqueness of the reduced contact form in Theorem \ref{thm:coiso-red}, since the characteristic leaves of $J^{-1}(\mu)$ are the orbits of the twisted action (Proposition \ref{prop:levelsets}(ii); recall that $G$ is connected). (iii) We have $H(g\star x)=H\big(\mathrm{Fl}^\Reeb_{\chi(g)}(\Phi_g(x))\big)$ and $H(\Phi_g(x))=H(x)$. As $\chi(G)=\R$ and $J^{-1}(\mu)$ is invariant under the flow of $\Reeb$ (Proposition \ref{prop:momentum}(iii)), $H|_{J^{-1}(\mu)}$ is invariant under the twisted action if and only if it is constant along the integral curves of $\Reeb$ in $J^{-1}(\mu)$, that is, if and only if $\Reeb(H)=0$ there. The expression of $\widetilde H_\mu$ follows from (ii). The statement about the dynamics is Proposition \ref{prop:coiso-dyn}, and $\Reeb(\widetilde H_\mu)\circ\psi_\mu=\Reeb(H)|_{J^{-1}(\mu)}=0$ because $\Reeb$ is $\psi_\mu$-related to the Reeb vector field of $\eta_{Q/G}$.
\end{proof}

\begin{remark}
Proposition \ref{prop:shift} is a contact version of the shifting trick: the twisted reduction at $\mu$ is conjugated by the strict contactomorphism $T_\varphi$ to the reduction at zero, and the function $\varphi$ is absorbed by the variable $z$, so that no magnetic term appears. In the symplectic case, for $\mu$ vanishing on $[\g,\g]$, the form $\langle\mu,\mathcal A\rangle-d\varphi$ is basic for any principal connection $\mathcal A$, so that the magnetic term $B_\mu$ of \cite{dLMdDV17} is exact and can also be removed by a fibre translation. The genuinely contact feature is the hypothesis $\Reeb(H)=0$ on $J^{-1}(\mu)$, without which $H$ does not descend to the reduced space.
\end{remark}

\begin{theorem}[Twisted reduction of the Hamilton--Jacobi problem]\label{thm:HJ-twisted}
Let $H$ be a $G$-invariant Hamiltonian with $\Reeb(H)=0$ on $J^{-1}(\mu)$, and let $\widetilde H_\mu$ be as in Proposition \ref{prop:shift}.
\begin{enumerate}[label=(\roman*)]
\item Every closed connected Legendrian submanifold $L$ with $J|_L=\mu$ is invariant under the twisted action, and $L\mapsto\psi_\mu(L)$ is a one-to-one correspondence between closed Legendrian solutions of the Hamilton--Jacobi problem for $\mathcal{E}_H$ (resp.\ $X_H$) contained in $J^{-1}(\mu)$ and closed Legendrian solutions of the Hamilton--Jacobi problem for $\mathcal{E}_{\widetilde H_\mu}$ (resp.\ $X_{\widetilde H_\mu}$), which preserves the value $E$ of the Hamiltonian.
\item $j^1S(Q)\subseteq J^{-1}(\mu)$ if and only if $S=\varphi+\widetilde S\circ\pi_G$ for a function $\widetilde S\colon Q/G\to\R$. In this case $\psi_\mu(j^1S(Q))=j^1\widetilde S(Q/G)$, and the reduced Hamilton--Jacobi equation is
\[
\widetilde H_\mu\big(j^1\widetilde S\big)=E\quad(\text{resp. }=0),\quad\text{i.e.,}\quad H\Big((T_q\pi_G)^*d\widetilde S+d\varphi(q),\,\widetilde S(\pi_G(q))+\varphi(q)\Big)=E .
\]
\item Let $G$ act on $Q\times\R$ by $g\cdot(q,z)=(gq,z+\chi(g))$ and let $\varpi_\mu(q,z)=(\pi_G(q),z-\varphi(q))$. Sections $\gamma$ of $\pi$ with coisotropic image, with values in $J^{-1}(\mu)$ and equivariant (i.e., $\gamma(g\cdot(q,z))=g\star\gamma(q,z)$), are in one-to-one correspondence with sections $\tilde\gamma$ of $T^*(Q/G)\times\R\to(Q/G)\times\R$ with coisotropic image, through $\tilde\gamma\circ\varpi_\mu=\psi_\mu\circ\gamma$; explicitly,
\[
\gamma(q,z)=\Big((T_q\pi_G)^*\tilde\gamma_{Q/G}\big(\pi_G(q),z-\varphi(q)\big)+d\varphi(q),\,z\Big).
\]
Moreover, $\gamma$ is a solution of the Hamilton--Jacobi problem for $X_H$ (resp.\ $\mathcal{E}_H$) if and only if $\tilde\gamma$ is a solution for $X_{\widetilde H_\mu}$ (resp.\ $\mathcal{E}_{\widetilde H_\mu}$).
\end{enumerate}
\end{theorem}

\begin{proof}
(i) By Theorem \ref{thm:invariant-legendrian}(ii), the twisted generators are tangent to $L$; since $L$ is closed and $G$ is connected, $L$ is invariant under the twisted action. By Proposition \ref{prop:shift}(i), $T_\varphi(L)$ is a closed $G$-invariant Legendrian submanifold contained in $J^{-1}(0)$. The proofs of Theorems \ref{thm:HJ-red}(ii) and \ref{thm:HJ-rec} only use that $\psi_0$ is a surjective submersion whose fibres are the $G$-orbits and that $\psi_0^*\eta_{Q/G}=i^*\eta_Q$; applied to $\psi_\mu=\psi_0\circ T_\varphi$, they show that $\psi_\mu(L)$ is a closed Legendrian submanifold and that, conversely, $\psi_\mu^{-1}(\widetilde L)$ is a closed Legendrian submanifold contained in $J^{-1}(\mu)$ for every closed Legendrian $\widetilde L$. Since $H=\widetilde H_\mu\circ\psi_\mu$ on $J^{-1}(\mu)$, $H|_L=E$ if and only if $\widetilde H_\mu|_{\psi_\mu(L)}=E$, and the assertion follows from Lemma \ref{lem:tangency}.
(ii) $J(j^1S(q))=\mu$ for all $q$ if and only if $dS(\xi_Q)=\mu(\xi)=d\varphi(\xi_Q)$ for all $\xi$, that is, if and only if $S-\varphi$ is $G$-invariant, i.e., $S-\varphi=\widetilde S\circ\pi_G$. Then $T_\varphi(j^1S(q))=\big(dS(q)-d\varphi(q),S(q)-\varphi(q)\big)=j^1(\widetilde S\circ\pi_G)(q)$, and $\psi_0$ maps it onto $j^1\widetilde S(\pi_G(q))$ (Theorem \ref{thm:HJ-red}(iii)).
(iii) $\varpi_\mu$ is a surjective submersion whose fibres are the orbits of the action on $Q\times\R$, $\tilde\pi\circ\psi_\mu=\varpi_\mu\circ\pi$ on $J^{-1}(\mu)$, where $\tilde\pi$ is the projection $T^*(Q/G)\times\R\to(Q/G)\times\R$, and $\partial_z$ is $\varpi_\mu$-related to $\partial_z$. The proof of Theorem \ref{thm:HJ-red-coiso} then applies verbatim, with $\psi_0$ and $\pi_G\times\id_\R$ replaced by $\psi_\mu$ and $\varpi_\mu$; in particular $\varpi_\mu^*\vartheta_{\tilde\gamma}=\gamma^*\psi_\mu^*\eta_{Q/G}=\vartheta_\gamma$ and $H\circ\gamma=\widetilde H_\mu\circ\tilde\gamma\circ\varpi_\mu$. The explicit expression of $\gamma$ is obtained by inverting $\psi_\mu$.
\end{proof}

\begin{remark}
By Corollary \ref{cor:no-nonzero}, $\Reeb(H)$ vanishes on every Legendrian solution with constant non-zero momentum, so the hypothesis of Theorem \ref{thm:HJ-twisted} is only slightly stronger than a necessary condition; it is needed for $H$ to descend to the reduced space. For a cyclic coordinate, Theorem \ref{thm:HJ-twisted}(ii) is exactly the classical ansatz $S=\beta q^1+\widetilde S(q')$ (Section \ref{sec:examples}).
\end{remark}

\subsection{Reduction along rays: the dissipative case}

Let $J_K\colon M\to\mathfrak k^*$ be the momentum map of the (lifted) action of $K$, i.e., $J_K(x)=J(x)|_{\mathfrak k}$, and consider the preimage of the ray $\R^+\mu=\{s\mu\mid s>0\}$,
\[
C_+=J^{-1}(\R^+\mu)=\{x\in M\mid J(x)=s\mu\ \text{for some } s>0\}.
\]
Let $\psi_0^K\colon J_K^{-1}(0)\to T^*(Q/K)\times\R$ and $H_K$ be the map and the reduced Hamiltonian of Proposition \ref{prop:cotangent-red} for the action of $K$, and let $\hat\zeta_{Q/K}$ be the infinitesimal generator of $[\hat\zeta]\in\g/\mathfrak k$ for the action of $G/K$ on $Q/K$, so that $T\pi_K(\hat\zeta_Q)=\hat\zeta_{Q/K}\circ\pi_K$.

\begin{proposition}\label{prop:ray}
\begin{enumerate}[label=(\roman*)]
\item $C_+=\{x\in J_K^{-1}(0)\mid J_{\hat\zeta}(x)>0\}$, and $J=J_{\hat\zeta}\,\mu$ on $C_+$. Hence $C_+$ is an open subset of $J_K^{-1}(0)$; it is a coisotropic submanifold of codimension $r-1$, tangent to $\Reeb$, whose characteristic leaves are the $K$-orbits.
\item For every $G$-invariant Hamiltonian $H$, with no condition on $\Reeb(H)$, the vector fields $X_H$ and $\mathcal{E}_H$ are tangent to $C_+$.
\item $\psi_0^K$ maps $C_+$ onto the open subset
\[
U_+=\{(\beta,z)\in T^*(Q/K)\times\R\mid \langle\beta,\hat\zeta_{Q/K}\rangle>0\},
\]
so that the reduced space $C_+/K$ is contactomorphic to $(U_+,\eta_{Q/K})$. Since $\g_\mu=\g$ and $K$ is the connected subgroup with Lie algebra $\ker\mu$, this is the contact reduced space at the ray $\R^+\mu$ in the sense of Willett \cite{Wi02}.
\item The group $G/K\cong\R$ acts on $T^*(Q/K)\times\R$ by the cotangent lift of its action on $Q/K$, leaving $U_+$ invariant, with momentum map $\bar J(\beta,z)=\langle\beta,\hat\zeta_{Q/K}\rangle$, and $\bar J\circ\psi_0^K=J_{\hat\zeta}$. If $H$ is $G$-invariant, $H_K$ is $G/K$-invariant, $X_H$ and $\mathcal{E}_H$ are $\psi^K_0$-related to $X_{H_K}$ and $\mathcal{E}_{H_K}$ on $C_+$, and $X_{H_K}(\bar J)=\mathcal{E}_{H_K}(\bar J)=-\Reeb(H_K)\,\bar J$. No Legendrian submanifold contained in $U_+$ is invariant under the action of $G/K$.
\end{enumerate}
\end{proposition}

\begin{proof}
(i) $J(x)\in\R\mu$ if and only if $J(x)$ vanishes on $\mathfrak k=\ker\mu$, i.e., $J_K(x)=0$, and in this case $J(x)=J_{\hat\zeta}(x)\mu$ because $\mu(\hat\zeta)=1$. Since $K$ acts freely on $Q$, $0$ is a regular value of $J_K$ (Proposition \ref{prop:regular}), and Proposition \ref{prop:levelsets}(iv), applied to $K$, shows that $J_K^{-1}(0)$ is coisotropic, of codimension $\dim K=r-1$, tangent to $\Reeb$, with characteristic leaves the $K$-orbits ($K$ is connected). These properties are local, so they hold for the open subset $C_+$.
(ii) By Proposition \ref{prop:dissipated}, $X_H(J_\xi)=\mathcal{E}_H(J_\xi)=-\Reeb(H)J_\xi$, which vanishes on $J_K^{-1}(0)$ for $\xi\in\mathfrak k$; hence $X_H$ and $\mathcal{E}_H$ are tangent to $J_K^{-1}(0)$, and to its open subset $C_+$. (Alternatively, along an integral curve the momentum is multiplied by a positive factor.)
(iii) By Proposition \ref{prop:cotangent-red} applied to $K$, $\psi_0^K$ is a surjective submersion with fibres the $K$-orbits and $(\psi_0^K)^*\eta_{Q/K}=i^*\eta_Q$. For $\alpha_q=(T_q\pi_K)^*\beta$ we have $J_{\hat\zeta}(\alpha_q,z)=\beta(T_q\pi_K(\hat\zeta_Q(q)))=\langle\beta,\hat\zeta_{Q/K}\rangle$, so $\psi_0^K(C_+)=U_+$.
(iv) The action of $G/K$ on $Q/K$ is free and proper (Lemma \ref{lem:character}), and its cotangent lift has momentum map $\bar J$, as in Example \ref{ex:cotangent-lift}; the equality $\bar J\circ\psi^K_0=J_{\hat\zeta}$ was shown in (iii), and $U_+=\bar J^{-1}(\R^+)$ is invariant. The invariance of $H_K$ follows from that of $H$, and the relation between the dynamics from Proposition \ref{prop:cotangent-red} applied to $K$. The dissipation of $\bar J$ is Proposition \ref{prop:dissipated} for the action of $G/K$. Finally, by Theorem \ref{thm:invariant-legendrian}(i), a $G/K$-invariant Legendrian submanifold is contained in $\bar J^{-1}(0)$, which does not meet $U_+$.
\end{proof}

Since $C_+$ is an open subset of $J_K^{-1}(0)$, the next theorem is essentially Theorems \ref{thm:HJ-red}, \ref{thm:HJ-rec} and \ref{thm:HJ-red-coiso} applied to the subgroup $K$ and restricted to the open sets $C_+$ and $U_+$. What is specific to the dissipative setting is that $C_+$ is invariant under the dynamics of every $G$-invariant Hamiltonian, with no condition on $\Reeb(H)$ (Proposition \ref{prop:ray}(ii)), that $C_+/K$ is Willett's reduced space at the ray, and that the residual symmetry $G/K$ cannot be used further (Remark \ref{rem:ray}(a)).

\begin{theorem}[Reduction of the Hamilton--Jacobi problem along rays]\label{thm:HJ-ray}
Let $H$ be a $G$-invariant Hamiltonian.
\begin{enumerate}[label=(\roman*)]
\item Every Legendrian submanifold $L$ which is closed in $C_+$ is $K$-invariant, and $L\mapsto\psi^K_0(L)$ is a one-to-one correspondence between Legendrian solutions of the Hamilton--Jacobi problem for $\mathcal{E}_H$ (resp.\ $X_H$) closed in $C_+$, i.e., on which the momentum map takes values in $\R^+\mu$, and Legendrian solutions for $\mathcal{E}_{H_K}$ (resp.\ $X_{H_K}$) closed in $U_+$, which preserves the value $E$ of the Hamiltonian.
\item For an open $K$-invariant subset $Q_0\subseteq Q$, $J(j^1S(Q_0))\subseteq\R^+\mu$ if and only if $S=S_K\circ\pi_K$ with $S_K\colon\pi_K(Q_0)\to\R$ such that $dS_K(\hat\zeta_{Q/K})>0$. The reduced Hamilton--Jacobi equation $H_K(j^1S_K)=E$ (resp.\ $=0$) is an equation on $Q/K$, which has $r-1$ fewer independent variables than $Q$.
\item The analogue of Theorem \ref{thm:HJ-red-coiso} holds for $K$-equivariant sections of $\pi$ with coisotropic image and values in $C_+$, which correspond to sections of $T^*(Q/K)\times\R\to(Q/K)\times\R$ with coisotropic image and values in $U_+$.
\end{enumerate}
\end{theorem}

\begin{proof}
(i) Since $L\subseteq J_K^{-1}(0)$, the infinitesimal generators of the $K$-action are tangent to $L$ (Theorem \ref{thm:invariant-legendrian}(i)), and, $L$ being closed in the $K$-invariant open set $C_+$ of $J_K^{-1}(0)$ and $K$ connected, $L$ is $K$-invariant. Now $H$ is $K$-invariant and $K$ acts freely and properly on $Q$, so Theorems \ref{thm:HJ-red} and \ref{thm:HJ-rec} for $K$ (whose proofs are local over the orbit space and apply to the invariant open subsets $C_+$ and $U_+$) give the correspondence; by Proposition \ref{prop:ray}(iii), $L\subseteq C_+$ if and only if $\psi^K_0(L)\subseteq U_+$. (ii) $J(j^1S(q))\in\R^+\mu$ if and only if $dS(\xi_Q)=0$ for $\xi\in\mathfrak k$ and $dS(\hat\zeta_Q)>0$, that is, if and only if $S$ is $K$-invariant ($K$ is connected) and $dS_K(\hat\zeta_{Q/K})>0$; the rest is Theorem \ref{thm:HJ-red}(iii) for $K$. (iii) is Theorem \ref{thm:HJ-red-coiso} applied to $K$, together with Proposition \ref{prop:ray}(iii).
\end{proof}

\begin{remark}\label{rem:ray}
\begin{enumerate}[label=(\alph*)]
\item At zero momentum the Hamilton--Jacobi equation loses $r$ variables (Theorem \ref{thm:HJ-red}); along a ray it only loses $r-1$, and the residual symmetry $G/K\cong\R$ cannot be used further, since no solution in $U_+$ is invariant under it (Proposition \ref{prop:ray}(iv)). This is the price of dissipation: the residual momentum $\bar J$ is not conserved but decays along the reduced dynamics.
\item If $r=1$, then $K$ is trivial and the reduction along rays is vacuous. This is the case of a single cyclic coordinate, where the time-dependent solutions $S(q,t)=\beta e^{-\lambda t}q^1+V(q',t)$ of Section \ref{sec:examples} have momentum on the ray $\R^+\beta$ but cannot be reduced.
\item The reduction along rays requires no assumption on $\Reeb(H)$. If $\Reeb(H)=0$ on $J^{-1}(\mu)\subseteq C_+$, the twisted reduction at $\mu$ can be recovered in two stages: first reduce along the ray by $K$, and then apply the twisted reduction to the residual action of $G/K\cong\R$ on $U_+$ at the value $\bar J=1$. We do not pursue this here.
\item The reduction along rays is also the natural framework for the slices of Proposition \ref{prop:slices}, whose momenta $\mu(z)$ move along a ray.
\item Assumption (A) enters Proposition \ref{prop:ray} and Theorem \ref{thm:HJ-ray} only through the closedness of $K$. If the connected subgroup $K$ with Lie algebra $\ker\mu$ is closed but $\mu$ does not integrate to a homomorphism into $\R$ (e.g., $G$ a torus and $\mu$ rational), then $G/K\cong S^1$, Proposition \ref{prop:ray}(i)--(iii) and Theorem \ref{thm:HJ-ray} hold verbatim, and (iv) holds with $S^1$ in place of $\R$. For irrational $\mu$ on a torus, $K$ is not closed and the method does not apply.
\end{enumerate}
\end{remark}

\section{Examples}\label{sec:examples}

\subsection{Cyclic coordinates}

Let $Q=\R\times Q'$ with coordinates $(q^1,q')$ and assume that $H(q,p,z)$ does not depend on $q^1$. Then $H$ is invariant under the lift of the action of $G=\R$ by translations in $q^1$, whose momentum map is $J=p_1$. By Proposition \ref{prop:dissipated}, $\dot p_1=-p_1\,\partial H/\partial z$. By Theorem \ref{thm:HJ-red}, the $\R$-invariant Legendrian solutions lie in $\{p_1=0\}$ and are of the form $j^1S$ with $S=S(q')$, $T^*(Q/G)\times\R=T^*Q'\times\R$, and the reduced equation is
\[
\widetilde H\left(q',\pd{S}{q'},S\right)=H\left(q',p_1=0,\pd{S}{q'},S\right)=E .
\]
This is the contact analogue of the classical ansatz for cyclic coordinates \cite[Section 5.3.1]{dLMdDV17}, but only at zero momentum. In fact, the classical ansatz $S=\beta q^1+V(q')$ with $\beta\neq0$ gives
\[
H\big(q',\beta,\partial_{q'}V,\beta q^1+V\big)=E,
\]
which depends on $q^1$ unless $\partial H/\partial z=0$. For instance, for $H=H_0(q',p)+\lambda z$ the left-hand side is $H_0(q',\beta,\partial_{q'}V)+\lambda V+\lambda\beta q^1$, which is never constant if $\lambda\beta\neq0$. This is a particular case of Corollary \ref{cor:no-nonzero}: since $\Reeb(H)=\lambda\neq0$, there are no Legendrian solutions, invariant or not, on which $p_1$ takes a constant non-zero value, the momentum $p_1$ decaying as $e^{-\lambda t}$. In the conservative case $\lambda=0$, the Legendrian submanifold $j^1(\beta q^1+V)(Q)$ is not $G$-invariant (in agreement with Theorem \ref{thm:invariant-legendrian}(i)), but it is invariant under the twisted action of Remark \ref{rem:mu-nonzero}, $(q^1,q',p,z)\mapsto(q^1+s,q',p,z+\beta s)$, and the classical ansatz corresponds to reduction at $\mu=\beta$ by means of this twisted action. Indeed, with $\chi(s)=\beta s$ and $\varphi=\beta q^1$, Theorem \ref{thm:HJ-twisted}(ii) gives the solutions $S=\beta q^1+\widetilde S(q')$, with reduced Hamiltonian $\widetilde H_\beta(q',p',z)=H(q',\beta,p',z)$. Instead, Case 1 admits non-invariant time-dependent solutions of \eqref{eq:evolHJ} of the form $S(q,t)=\beta e^{-\lambda t}q^1+V(q',t)$, which follow the decay of the momentum: the terms in $q^1$ cancel and $V$ must satisfy $\partial_tV+H_0(q',\beta e^{-\lambda t},\partial_{q'}V)+\lambda V=0$. Their momenta lie on the ray $\R^+\beta$, but since $G=\R$ the reduction along rays is vacuous (Remark \ref{rem:ray}(b)); Section \ref{ex:ray} shows a case where it is not.

\subsection{A central force problem with linear dissipation}

Let $Q=\R^2\setminus\{0\}$ and consider the Hamiltonian on $T^*Q\times\R$
\[
H(q,p,z)=\frac12|p|^2+V(|q|)+\lambda z,\qquad\lambda>0,
\]
which describes a particle in a central potential subject to a friction force proportional to its velocity: the contact Hamilton equations give $\ddot q=-V'(|q|)\,q/|q|-\lambda\dot q$. The group $G=SO(2)$ acts freely on $Q$ by rotations, $H$ is invariant under the lifted action and the momentum map is the angular momentum $J=q^1p_2-q^2p_1$. By Proposition \ref{prop:dissipated}, $J$ decays as $e^{-\lambda t}$, and $J^{-1}(0)$ (radial motions) is invariant.

The quotient is $Q/G=(0,+\infty)$ with coordinate $r=|q|$, $\pi_G(q)=|q|$. A covector $\alpha_q\in J^{-1}(0)$ is of the form $\alpha_q=p_r\,dr$, and, since $|p|^2=p_r^2+J^2/r^2$, the reduced Hamiltonian on $T^*(0,\infty)\times\R$ is
\[
\widetilde H(r,p_r,z)=\frac12p_r^2+V(r)+\lambda z .
\]
By Theorems \ref{thm:HJ-red} and \ref{thm:HJ-rec}, the $SO(2)$-invariant horizontal Legendrian solutions for $\mathcal{E}_H$ are $j^1S$ with $S(q)=\widetilde S(|q|)$, where $\widetilde S$ satisfies the ordinary differential equation
\[
\frac12\widetilde S'(r)^2+V(r)+\lambda\widetilde S(r)=E,
\]
while the original Hamilton--Jacobi equation is the PDE $\frac12|\nabla S|^2+V(|q|)+\lambda S=E$ in two variables.

For the harmonic potential $V(r)=\frac{\kappa}{2}r^2$, we look for $\widetilde S(r)=a r^2+b$. The reduced equation becomes $\big(2a^2+\lambda a+\frac{\kappa}{2}\big)r^2+\lambda b=E$, so that
\[
a=\frac{-\lambda\pm\sqrt{\lambda^2-4\kappa}}{4},\qquad b=\frac{E}{\lambda},
\]
real in the overdamped case $\lambda^2\geq4\kappa$. The reconstructed solution is $S(q)=a|q|^2+E/\lambda$, with $L=j^1S(Q)=\{p=2a q,\ z=a|q|^2+E/\lambda\}$. The projected vector field is $Y(q)=\partial H/\partial p\circ j^1S=2aq$, whose integral curves $q(t)=e^{2at}q_0$ are exactly the radial exponential modes of the damped oscillator $\ddot q+\lambda\dot q+\kappa q=0$, whose characteristic roots are $2a$. Accordingly, $\mathcal{A}(Y)=0$ for the mechanical connection, so that $g(t)=e$ in Proposition \ref{prop:reconstruction} (see Remark \ref{rem:mech-conn}).

The same functions give solutions of the other cases. The family $S(q,c)=a|q|^2+c$ defines the coisotropic section $\gamma(q,z)=(q,2aq,z)$, and $H(j^1S_c)=\lambda c$; hence $\gamma$ is a solution of the Hamilton--Jacobi problem for $\mathcal{E}_H$ with $k(c)=\lambda c$ (Theorem \ref{thm:case3}) and also for $X_H$, with $k(c)=\lambda c$ in \eqref{eq:HJ1-leaves}. Reparametrizing as in Corollary \ref{cor:evolution} (on the region $c\neq0$, where $H\circ\gamma\neq0$), $dc/dt=-\lambda c$, we obtain the solution
\[
S(q,t)=a|q|^2+c_0e^{-\lambda t}
\]
of the evolutionary equation $\partial_tS+\frac12|\partial_qS|^2+\frac{\kappa}{2}|q|^2+\lambda S=0$ (with $c_0\neq0$), as can be checked directly. Finally, $E=0$ gives the Legendrian solution $S=a|q|^2$ for $X_H$. All these solutions are $SO(2)$-invariant and reduce to the corresponding solutions on $(0,\infty)$.

\subsection{Two particles on a line with linear dissipation}

Consider two particles on a line interacting through a potential depending on their relative position, with linear friction:
\[
H(q_1,q_2,p_1,p_2,z)=\frac12(p_1^2+p_2^2)+V(q_1-q_2)+\lambda z .
\]
For $\lambda=0$ and $V(x)=1/x^2$ this is the two-particle Calogero--Moser system studied in \cite[Section 5.4]{dLMdDV17}. The group $\R$ acts by simultaneous translations, $(q_1,q_2)\mapsto(q_1+s,q_2+s)$, with momentum map $J=p_1+p_2$ (the total momentum), which decays as $e^{-\lambda t}$. The quotient $Q/G=\R$ has coordinate $x=q_1-q_2$; a covector in $J^{-1}(0)$ is $\alpha=p\,dx=p\,dq_1-p\,dq_2$, so that $p_1=-p_2=p$ and the reduced Hamiltonian is
\[
\widetilde H(x,p,z)=p^2+V(x)+\lambda z .
\]
The $\R$-invariant solutions are $S(q_1,q_2)=\widetilde S(q_1-q_2)$ with
\[
\widetilde S'(x)^2+V(x)+\lambda\widetilde S(x)=E .
\]
For $\lambda=0$ the reduced equation becomes the conservative one, $\widetilde S'(x)^2+V(x)=E$, that is, the reduced Hamilton--Jacobi equation at zero total momentum of the conservative problem considered in \cite[Section 5.4]{dLMdDV17}. Notice that, by Corollary \ref{cor:no-nonzero}, for $\lambda\neq0$ there are no Legendrian solutions on which the total momentum takes a constant non-zero value. In the conservative case such solutions exist; they correspond in \cite{dLMdDV17} to the reduction at $J^{-1}(\mu)$, $\mu\neq0$, and in the contact setting to the twisted reduction of Theorem \ref{thm:HJ-twisted}. For $V(x)=\frac{\kappa}{2}x^2$ and $\lambda\neq0$, the solutions $\widetilde S=a x^2+E/\lambda$ with $4a^2+\lambda a+\frac{\kappa}{2}=0$ (real when $\lambda^2\geq8\kappa$) reconstruct to $S=a(q_1-q_2)^2+E/\lambda$. On $j^1S(Q)$ we have $p_1=-p_2=2a(q_1-q_2)$, hence $\dot x=4ax$, and $4a$ is a characteristic root of the equation $\ddot x+\lambda\dot x+2\kappa x=0$ of the relative motion.

\subsection{A particle in the plane with friction: reduction along rays}\label{ex:ray}

Let $Q=\R^2$ and $H(q,p,z)=\frac12|p|^2+\lambda z$, $\lambda>0$, which describes a free particle subject to linear friction, $\ddot q=-\lambda\dot q$. The group $G=\R^2$ acts by translations, $J=p$, and every $\mu\neq0$ vanishes on $[\g,\g]=0$. Let $\mu=\mathbf u$ be a unit vector, so that $\chi(s)=\mathbf u\cdot s$, $K=\{s\in\R^2\mid \mathbf u\cdot s=0\}$ and $\hat\zeta=\mathbf u$. Here $Q/K=\R$ with coordinate $w=\mathbf u\cdot q$, $\hat\zeta_{Q/K}=\partial_w$, and a covector in $J_K^{-1}(0)$ is of the form $p=p_w\,\mathbf u$. Hence $U_+=\{p_w>0\}\subseteq T^*\R\times\R$ and
\[
H_K(w,p_w,z)=\frac12p_w^2+\lambda z .
\]
Since $\Reeb(H)=\lambda\neq0$, there are no Legendrian solutions with constant non-zero momentum (Corollary \ref{cor:no-nonzero}), so neither Theorem \ref{thm:HJ-red} nor Theorem \ref{thm:HJ-twisted} gives solutions with $p\neq0$. Instead, by Theorem \ref{thm:HJ-ray}, the Legendrian solutions for $\mathcal{E}_H$ with momentum on the ray $\R^+\mathbf u$ are $j^1S$ with $S(q)=S_K(\mathbf u\cdot q)$, where
\[
\frac12S_K'(w)^2+\lambda S_K(w)=E,\qquad S_K'(w)>0 .
\]
The solutions are $S_K(w)=\frac{E}{\lambda}-\frac{\lambda}{2}(w-w_0)^2$ on $w<w_0$. Hence
\[
S(q)=\frac{E}{\lambda}-\frac{\lambda}{2}\,(\mathbf u\cdot q-w_0)^2,\qquad p=\partial_qS=-\lambda\,(\mathbf u\cdot q-w_0)\,\mathbf u\in\R^+\mathbf u,
\]
on the $K$-invariant half-plane $\{\mathbf u\cdot q<w_0\}$ (for $X_H$, take $E=0$). The projected vector field is $Y=p$, so that $w(t)=w_0+(w(0)-w_0)e^{-\lambda t}$ and $p(t)=p(0)e^{-\lambda t}$. The particle moves in the direction $\mathbf u$ and approaches the line $\mathbf u\cdot q=w_0$, at distance $|p(0)|/\lambda$ from its initial position, which is the stopping distance of the linearly damped motion. The momentum decays along the ray, as predicted by Propositions \ref{prop:dissipated} and \ref{prop:ray}(iv), and the reduced equation is an ordinary differential equation, whereas the original one is a partial differential equation in two variables. These solutions are neither $G$-invariant nor invariant under a twisted action, and they are the typical solutions with non-zero momentum of a dissipative system with symmetry. For $\lambda=0$ one recovers, by Theorem \ref{thm:HJ-twisted}, the solutions $S=|p_0|\,\mathbf u\cdot q+c$ of the free particle.

\section{Conclusions and future work}\label{sec:conclusions}

In this paper we have studied the interplay between the Hamilton--Jacobi theory of contact Hamiltonian systems, symmetries and coisotropic reduction, extending to the contact (and hence dissipative) setting the results of \cite{dLMdDV17} for symplectic Hamiltonian systems.

First, we have clarified the geometry of the four contact Hamilton--Jacobi equations of \cite{dLLMB21}. Sections of $T^*Q\times\R\to Q\times\R$ with coisotropic image are exactly the Legendrian foliations, i.e., one-parameter families of $1$-jets $j^1S_c$, so that the assumption of Lagrangian slices is not necessary. With this observation, the Hamilton--Jacobi equation for the evolution vector field says that the Hamiltonian is constant on each Legendrian leaf, and the one for the Hamiltonian vector field is equivalent to the evolutionary contact Hamilton--Jacobi equation $\partial_tS+H(q,\partial_qS,S)=0$, the leaves being transported by the contact flow. The equations for Legendrian sections admit an intrinsic characterization in terms of the tangency of $X_H$ or $\mathcal{E}_H$ to Legendrian submanifolds, which allows us to define generalized (non-horizontal) solutions.

Second, we have analysed momentum maps of actions by contactomorphisms. They are dissipated quantities (so that, in the genuinely dissipative case $\Reeb(H)\neq0$, no Legendrian solution has constant non-zero momentum), and a closed connected Legendrian submanifold is invariant if and only if the momentum map vanishes on it; constant non-zero values correspond to invariance under a twisted action involving the Reeb flow. Using the contact coisotropic reduction theorem for coisotropic submanifolds tangent to the Reeb vector field and the reduction of Legendrian submanifolds \cite{dLIL24}, we have recovered the contact reduction at zero momentum and the reduction of the contact dynamics.

Third, we have obtained a reduction and reconstruction procedure for the four contact Hamilton--Jacobi equations: $G$-invariant solutions for a $G$-invariant Hamiltonian on $T^*Q\times\R$ are in one-to-one correspondence with solutions for the reduced Hamiltonian on $T^*(Q/G)\times\R$. A notable difference with the symplectic case is that $G$-invariant solutions only occur at the zero level of the momentum map, no magnetic terms appear, and the generating functions of invariant solutions are invariant. Non-zero momenta $\mu$ vanishing on $[\g,\g]$ are treated in Section \ref{sec:nonzero}. In the conservative case $\Reeb(H)=0$, a strict contactomorphism conjugates the twisted reduction at $\mu$ with the reduction at zero, and the classical ansatz for cyclic variables is recovered. In the dissipative case, the preimage of the ray $\R^+\mu$ is invariant and is reduced by the subgroup with Lie algebra $\ker\mu$ to Willett's reduced space; the reduced Hamilton--Jacobi equation loses one variable fewer than at zero momentum, which is the price of dissipation. The examples show that, in the presence of dissipation, the classical ansatz for cyclic variables only survives at zero momentum, while the reduction along rays captures the solutions with decaying momentum.

Several problems remain open:
\begin{enumerate}[label=(\arabic*)]
\item \emph{General symmetries.} We have considered cotangent lifted actions, which act trivially on $z$. It would be interesting to extend the theory to actions by conformal contactomorphisms (in particular, scaling symmetries \cite{BJS23}) and to non-lifted actions, as well as to singular momentum values, using the theory of singular contact reduction \cite{Wi02,ZZ06}.
\item \emph{General momentum values.} Our reduction at non-zero momentum assumes that $\mu$ vanishes on $[\g,\g]$ and integrates to a character of $G$. For general $\mu$, $J^{-1}(\mu)$ is not coisotropic and Willett's reduction uses the subgroup with Lie algebra $\ker\mu\cap\g_\mu$; extending Theorems \ref{thm:HJ-twisted} and \ref{thm:HJ-ray} to this setting, where magnetic terms may appear, and making precise the reduction by stages of Remark \ref{rem:ray}(c), deserve further study.
\item \emph{Complete solutions and integrability.} The reduction of complete solutions and its relation with the contact Liouville--Arnold theorem \cite{CdLLLG23} and with the integrability results of \cite{dLLLGR23,GP20} and with the contact bi-Hamiltonian structures of \cite{CdLEILG25}.
\item \emph{Time-dependent systems.} The extension to cocontact manifolds, describing time-dependent dissipative systems, and to the evolutionary equation in the sense of viscosity solutions \cite{CCJWY20}.
\item \emph{Numerical methods.} The construction of contact geometric integrators from (reduced) generating functions, in the spirit of the discrete Herglotz principle, contact variational integrators \cite{VBS19,BSVZ20} and the discrete contact Hamilton--Jacobi theory \cite{ESZ24}, and its application to thermodynamic systems \cite{SdLLMdD20,dLB25}.
\item \emph{Field theories.} The extension of these results to multicontact field theories, and to $k$-contact field theories, for which Hamilton--Jacobi theory \cite{dLLRS26} and Marsden--Meyer--Weinstein reduction \cite{dLRVZ25} have recently been developed.
\end{enumerate}

\section*{Acknowledgements}
Manuel de Le\'on acknowledges the financial support of the Ministerio de Ciencia, Innovaci\'on y Universidades (Spain), grant PID2022-137909NB-C21, and of the Severo Ochoa Programme for Centers of Excellence in R\&D (CEX2023-001347-S). He also acknowledges the partial support of the National Natural Science Foundation of China (Grant No.\ W2541005), of the Guangdong Provincial Key Laboratory of Mathematical and Neural Dynamical Systems (Grant No.\ 2024B1212010004), and of a Distinguished Visiting Professorship at Great Bay University.
Mar\'ia Emma Eyrea Iraz\'u acknowledges the support of the Universidad Nacional de La Plata.

\end{document}